\documentclass[10pt]{article}
\usepackage[dvips]{graphicx}
\usepackage{amsmath,amsfonts,amssymb,latexsym,epsfig}
\usepackage{mathrsfs}
\usepackage{verbatim}
\usepackage{latexsym}
\usepackage{amsthm}
\usepackage{amssymb}
\usepackage{graphics}
\usepackage{amsbsy}
\usepackage{times}

\newtheorem{theorem}{Theorem}[section]
\newtheorem{lemma}[theorem]{Lemma}
\newtheorem{remark}[theorem]{Remark}

\newtheorem{prop}[theorem]{Proposition}
\newtheorem{assumptions}[theorem]{Assumptions}
\newtheorem{assumption}[theorem]{Assumption}

\numberwithin{equation}{section}
\graphicspath{{../figures/}}

\newcommand{\E}{{\mathbb E}}
\newcommand{\bbE}{{\mathbb E}}
\newcommand{\bbP}{{\mathbb P}}

\newcommand{\T}{{\mathbb T}}
\newcommand{\R}{{\mathbb R  }}
\newcommand{\bL}{\mathbb L}

\newcommand{\oT}{\hat{\Theta}}

\newcommand{\bbT}{\mathbb{T}}
\newcommand{\bbR}{\mathbb{R}}

\newcommand{\cL}{\mathcal L}

\newcommand{\cX}{\mathcal X}
\newcommand{\cY}{\mathcal Y}
\newcommand{\cZ}{\mathcal Z}

\newcommand{\eps}{\epsilon}

\newcommand{\epsr}{\sqrt{\epsilon}}
\newcommand{\epss}{\eps^2}
\newcommand{\rato}{\frac{1}{\eps}}
\newcommand{\ratt}{\frac{1}{\eps^2}}
\newcommand{\rath}{\frac{1}{\epsr}}
\newcommand{\htheta}{\hat{\theta}}

\newcommand{\Phie}{\psi^{\epsilon}}
\newcommand{\pie}{\rho^{\epsilon}}

\begin{document}
\setlength{\baselineskip}{10pt}
\title{MAXIMUM LIKELIHOOD DRIFT  ESTIMATION FOR MULTISCALE DIFFUSIONS}
\author{
A.Papavasiliou
\footnote{E-mail address: a.papavasiliou@warwick.ac.uk.} \\
        Statistics Department\\
        Warwick University \\
        Coventry CV4 7AL, UK\\
\\
G.A. Pavliotis\footnote{
E-mail address: g.paviotis@imperial.ac.uk.} \\
        Department of Mathematics\\
    Imperial College London \\
        London SW7 2AZ, UK \\
        and \\
        A.M. Stuart
\footnote{E-mail address: a.m.stuart@warwick.ac.uk.} \\
        Mathematics Institute \\
        Warwick University \\
        Coventry CV4 7AL, UK
                    }
\maketitle

\begin{abstract}
We study the problem of parameter estimation using maximum likelihood for fast/slow
systems of stochastic differential equations. Our aim is to shed light on the problem of
model/data mismatch at small scales. We consider two classes of fast/slow problems for
which a closed coarse-grained equation for the slow variables can be rigorously derived,
which we refer to as averaging and homogenization problems. We ask whether, given data
from the slow variable in the fast/slow system, we can correctly estimate parameters in
the drift of the coarse-grained equation for the slow variable, using maximum likelihood.
We show that, whereas the maximum likelihood estimator is asymptotically unbiased for the
averaging problem, for the homogenization problem maximum likelihood fails unless we
subsample the data at an appropriate rate. An explicit formula for the asymptotic error
in the log likelihood function is presented. Our theory is applied to two simple examples
from molecular dynamics.
\end{abstract}


\section{Introduction}

Fitting stochastic differential equations (SDEs) to time-series data is often a useful
way of extracting simple model fits which capture important aspects of the
dynamics~\cite{GKS04}. However, whilst the data may well be compatible with an SDE model
in many respects, it is often incompatible with the desired model at small scales. Since
many commonly applied statistical techniques see the data at small scales this can lead
to inconsistencies between the data and the desired model fit. This phenomenon appears
quite often in econometrics~\cite{AitMykZha05b, AitMykZha05a, OlhPavlSyk08}, where the
term market microstructure noise is used to describe the high frequency/small scale part
of the data as well as in molecular dynamics~\cite{PavlSt06}. In essence, the problem
that we are facing is that there is an inconsistency between the coarse-grained model
that we are using and the microscopic dynamics from which the data is generated, at small
scales. Similar problems appear quite often in statistical inference, in the context of
parameter estimation for misspecified or incorrect models~\cite[Sec. 2.6]{Kut04}.

The aim of this paper is to create a theoretical framework in which it is possible to
study this issue, in order to gain better insight into how it is manifest in practice,
and how to overcome it. In particular our goal is to investigate the following problem:
how can we fit data obtained from the high-dimensional, multiscale full dynamics to a
low-dimensional, coarse grained model which governs the evolution of the resolved
("slow") degrees of freedom? We will study this question for a class of stochastic
systems for which we can derive rigorously a coarse grained description for the dynamics
of the resolved variables. More specifically,  we will work in the framework of coupled
systems of multiscale SDEs for a pair of unknown functions $(x(t),y(t))$. We assume that
$y(t)$ is fast, relative to $x(t)$, and that the equations average or homogenize to give
a closed equation for $X(t)$ to which $x(t)$ converges in the limit of infinite scale
separation. The function $X(t)$ then approximates $x(t)$, typically in the sense of weak
convergence of probability measures~\cite{EthKur86, PavlSt08}. We then ask the following
question: given data for $x(t)$, from the coupled system, can we correctly identify
parameters in the averaged or homogenized model for $X(t)$?

Fast/slow systems of SDEs of this form have been studied extensively over the last four
decades~\cite{lions, PapStrVar77, PavlSt08} and the references therein. Recently, various
methods have been proposed for solving numerically these
SDEs~\cite{ELV05,GivKevKup06,Vand03}. In these works, the coefficients of the limiting
SDE are calculated "on the fly" from simulations of the fast/slow system. There is a
direct link between these numerical methods and our approach in that our goal is also to
infer information about the coefficients in the coarse-grained equation using data from
the multiscale system. However, our interest is mainly in situations where the
"microscopic" multiscale system is not known explicitly. From this point of view, we
merely use the multiscale stochastic system as our "data generating process"; our goal is
to fit this data to the coarse-grained equation for $X(t)$, the limit of the slow
variable $x(t)$.

A first step towards the understanding of this problem was taken in~\cite{PavlSt06}.
There, the data generating process $x(t)$ was taken to be the path of a particle moving
in a multiscale potential under the influence of thermal noise. The goal was to identify
parameters in the drift as well as the diffusion coefficient  in the homogenized model
for $X(t)$, the weak limit of $x(t)$. It was shown that the maximum likelihood estimator
is asymptotically biased and that subsampling is necessary in order to estimate the
parameters of the homogenized limit correctly, based on a time series (i.e. single
observation) of $x(t)$.

In this paper we extend the analysis to more general classes of fast/slow
systems of SDEs for which either an averaging or homogenization principle
holds~\cite{PavlSt08}. We consider cases where the drift in the averaged or
homogenized equation contains parameters which we want to estimate using
observations of the slow variable in the fast/slow system.
 We show that in the case of averaging the maximum likelihood function is
asymptotically unbiased and that we can estimate correctly the parameters of the drift in the averaged model from a single path of the slow variable
$x(t)$. On the other hand, we show
rigorously that the maximum likelihood estimator is
asymptotically biased for homogenization problems. In particular, an additional term
appears in the likelihood function in the limit of infinite scale separation.
We show then that this term vanishes, and hence that the maximum
likelihood estimator becomes asymptotically unbiased, provided that we subsample at an appropriate
rate.

To be more specific, in this paper we will consider fast/slow systems of SDEs of the form
\begin{subequations}
\label{e:av_intro}
\begin{eqnarray}\label{e:av_intro_x}
\frac{dx}{dt} &=& f_1(x,y)+\alpha_0(x,y)\frac{dU}{dt}+\alpha_1(x,y)\frac{dV}{dt}, \\
\frac{dy}{dt} &=& \rato g_0(x,y)+\rath \beta(x,y)\frac{dV}{dt};
\end{eqnarray}
\end{subequations}
or the SDEs
\begin{subequations}
\label{e:hom_intro}
\begin{eqnarray}\label{e:hom_intro_x}
\frac{dx}{dt} &=& \rato
f_0(x,y)+f_1(x,y)+\alpha_0(x,y)\frac{dU}{dt}+\alpha_1(x,y)\frac{dV}{dt}, \\
\frac{dy}{dt} & =& \ratt g_0(x,y)+\rato g_1(x,y)+\rato \beta(x,y)\frac{dV}{dt}.
\end{eqnarray}
\end{subequations}
We will refer to equations~\eqref{e:av_intro} as the {\bf averaging problem} and to
equations~\eqref{e:hom_intro} as the {\bf homogenization problem}. In both cases our
assumptions on the coefficients in the SDEs are such that a coarse-grained (averaged or
homogenized) equation exists, which is of the form
\begin{equation}\label{e:coarse_intro}
\frac{d X}{d t} = F(X; \theta) + K(X) \frac{d W}{d t}.
\end{equation}
The slow variable $x(t)$ converges weakly, in the limit as $\eps \rightarrow 0$, to
$X(t)$, the solution of~\eqref{e:coarse_intro}. We assume that the vector field
$F(X;\theta)$ depends on a set of parameters  $\theta$ that we want to estimate based on data from
either the averaging or the homogenization problem. We suppose that the actual
drift compatible with the data is given by $F(X)=F(X;\theta_0).$ We ask whether it is
possible to  correctly identify $\theta=\theta_0$ by finding the {\em maximum likelihood
estimator} (MLE) when using a statistical model of the form \eqref{e:coarse_intro}, but
given data from \eqref{e:av_intro} or \eqref{e:hom_intro}. Our main results can be
stated, informally, as follows.
\begin{theorem}
Assume that we are given continuous time data. The MLE for the averaging problem (i.e.
fitting data from~\eqref{e:av_intro_x} to~\eqref{e:coarse_intro}) is asymptotically
unbiased. On the other hand, the MLE for the homogenization problem (i.e. fitting data
from~\eqref{e:hom_intro_x} to~\eqref{e:coarse_intro}) is asymptotically biased and an
explicit formula for the asymptotic error in the likelihood, $E_{\infty}$, can be
obtained.
\end{theorem}
Precise statements of the above results can be found in Theorems~\ref{th:avpe},
\ref{thm:mle_av} and \ref{th:hompe}.

The failure of the MLE when applied to the homogenization problem is due to the presence
of high frequency data. Naturally, in order to be able to identify correctly the
parameter $\theta = \theta_0$ in~\eqref{e:coarse_intro} using data
from~\eqref{e:hom_intro_x} subsampling at an appropriate rate is necessary.
\begin{theorem}
The MLE for the homogenization problem becomes asymptotically unbiased if we subsample at
an appropriate rate.
\end{theorem}
Roughly speaking, the sampling rate should be between the two characteristic time scales
of the fast/slow SDEs~\eqref{e:hom_intro}, $1$ and $\eps^2$. The precise statement of
this result can be found in Theorems~\ref{th:main1} and~\ref{thm:subsam}.
IIn practice real data will not come explicitly from a scale-separated
model like \eqref{e:av_intro_x} or \eqref{e:hom_intro_x}. However real data
is often multiscale in character. Thus the results in this paper shed light
on the pitfalls that may arise when fitting simplified statistical models
to multiscale data. Furthermore the results indicate the central, and
subtle, role played by subsampling data in order to overcome mismatch
between model and data at small scales.

The rest of the paper is organized as follows. In Section~\ref{sec:set_up}
we study the fast/slow
stochastic systems introduced above, and prove appropriate averaging and
homogenization theorems. In Section~\ref{sec:par_est} we introduce
the maximum likelihood function for \eqref{e:coarse_intro}
and study its limiting behavior, given
data from the averaging and
homogenization problems \eqref{e:av_intro_x} and \eqref{e:hom_intro_x}.
In Section~\ref{sec:subsamp} we show that, when subsampling at
an appropriate rate, the maximum likelihood estimator for the homogenization problem
becomes asymptotically unbiased. In Section~\ref{sec:examp} we present examples of
fast/slow stochastic systems that fit into the general framework of this paper.
Section~\ref{sec:conc} is reserved for conclusions. Various technical results are proved
in the appendices.

\section{Set-Up}\label{sec:set_up}

We will consider fast/slow systems of SDEs for the variables $(x, \, y) \in
{\cal X} \times {\cal Y}$. We can take, for example, ${\cal X}\times {\cal
Y} = \R^l \times \R^{d-l}$ or ${\cal X}\times {\cal
Y} = \T^l \times \T^{d-l}$. In the second case, where the state space is compact,
all of the assumptions that we need for the proofs of our results can be
justified using elliptic PDEs theory.

Let $\varphi_{\xi}^t(y)$ denote the Markov process which solves the SDE
\begin{equation}
\label{eq:fast}
\frac{d}{dt}\Bigl(\varphi_{\xi}^t(y)\Bigr)=g_0\bigl(\xi,\varphi_{\xi}^t(y)\bigr)
+\beta\bigl(\xi,\varphi^t_{\xi}(y)\bigr)\frac{dV}{dt}, \quad \varphi^0_{\xi}(y)=y.
\end{equation}
Here $\xi \in \cX$ is a fixed parameter and, for each $t \ge 0$, $\varphi_{\xi}^t(y) \in
\cY$, $g_0:\cX \times \cY \to \bbR^{d-l}$, $\beta: \cX \times \cY \to \bbR^{(d-l) \times
m}$ and $V$ is a standard Brownian motion in $m$ dimensions.\footnote{Throughout this
paper we write stochastic differential equations as identities in fully differentiated
form, even though Brownian motion is not differentiable. In all cases the identity should
be interpeted as holding in integrated form, with the It\^{o} interpreation of the
stochastic integral.} The generator of the process is
\begin{equation}
{\cal L}_0(\xi)=g_0(\xi,y) \cdot \nabla_y+\frac12 B(\xi,y):\nabla_y \nabla_y
\label{eq:ygen}
\end{equation}
with $B(\xi,y):=\beta(\xi,y)\beta(\xi,y)^T.$
Notice that $\cL_0(\xi)$ is a differential operator in $y$ alone,
with $\xi$ a parameter.

Our interest is in data generated by the projection onto the
$x$ coordinate of  systems of SDEs  for $(x,y)$ in ${\cal X}
\times {\cal Y}.$
In particular, for $U$ a standard Brownian motion in $\bbR^n$ we
will consider either of the following coupled systems of SDEs:
\begin{subequations}
\label{eq:av}
\begin{eqnarray}
\frac{dx}{dt} &=& f_1(x,y)+\alpha_0(x,y)\frac{dU}{dt}+\alpha_1(x,y)\frac{dV}{dt}, \\
\frac{dy}{dt} &=& \rato g_0(x,y)+\rath \beta(x,y)\frac{dV}{dt};
\end{eqnarray}
\end{subequations}
or the SDEs
\begin{subequations}
\label{eq:hom}
\begin{eqnarray}
\frac{dx}{dt} &=& \rato
f_0(x,y)+f_1(x,y)+\alpha_0(x,y)\frac{dU}{dt}+\alpha_1(x,y)\frac{dV}{dt}, \\
\frac{dy}{dt} & =& \ratt g_0(x,y)+\rato g_1(x,y)+\rato \beta(x,y)\frac{dV}{dt}.
\end{eqnarray}
\end{subequations}
Here $f_i: \cX \times \cY \to \bbR^l, \alpha_0: \cX \times \cY \to
\bbR^{l\times n}, \alpha_1: \cX \times \cY \to \bbR^{l \times m},$
$g_1: \cX \times \cY \to \bbR^{d-l}$ and $g_0, \beta$ and $V$ are as above.

\begin{assumptions}
\label{ass:a1}

\begin{itemize}

\item The equation
$$-\cL_0^*(\xi) \rho(y;\xi)=0, \quad \int_{\cY} \rho(y;\xi)dy=1$$
has a unique non-negative solution $\rho(y;\xi) \in L^1(\cY)$
for every $\xi \in \cX$; furthermore $\rho(y;\xi)$ is $C^{\infty}$
in $y$ and $\xi.$

\item For each $\xi \in {\cal X}$ define the weighted Hilbert
space $L^2_{\rho}(\cY;\xi)$ with inner-product
$$\langle a,b \rangle_{\rho}:=\int_{\cY} \rho(y;\xi)a(y)b(y)dy.$$
For all $\xi \in {\cal X}$ Poisson equation
$$-\cL_0(\xi) \Theta(y;\xi)=h(y;\xi), \quad \int_{\cY} \rho(y;\xi) \Theta(y;\xi)dy=0$$
has a unique  solution $\Theta(y;\xi) \in L^2_{\rho}(\cY;\xi),$
provided that
$$\int_{\cY}\rho(y;\xi)h(y;\xi)dy=0.$$

\item The functions $f_i, g_i, \alpha_i, \beta$ and all derivatives
are uniformly bounded in $\cX \times \cY.$

\item If $h(y;\xi)$ and all its derivatives with respect to $y,\xi$
are uniformly bounded in $\cX \times \cY$
then the same is true of $\Theta$ solving the Poisson equation above.

\end{itemize}

\end{assumptions}

\begin{remark}
In the case where the state space of the fast process is compact, $\mathcal{Y} = \T^{d
-\ell}$, and the diffusion matrix $B(\xi, y)$ is positive definite the above assumptions
can be easily proved using elliptic PDE theory~\cite[Ch. 6]{PavlSt08}. Similar results
can also be proved without the compactness and uniform ellipticity
assumptions~\cite{PV01, PV03, PV05}.
\end{remark}

The first assumption essentially states the the process \eqref{eq:fast}
is ergodic, for each $\xi \in {\cal X}.$ Let $\cL_0=\cL_0(x)$ and define
\begin{align*}
\cL_1&=f_0 \cdot \nabla_x+g_1 \cdot \nabla_y+C:\nabla_y \nabla_x,\\
\cL_2&=f_1 \cdot \nabla_x+\frac12 A: \nabla_x \nabla_x,
\end{align*}
where
\begin{align*}
A(x,y)&=\alpha_0(x,y)\alpha_0(x,y)^T+\alpha_1(x,y)\alpha_1(x,y)^T,\\
C(x,y)&=\alpha_1(x,y)\beta(x,y)^T.
\end{align*}
The generators for the Markov processes defined by equations \eqref{eq:av}
and \eqref{eq:hom} respectively are
\begin{align}
{\cal L}_{av}&=\frac{1}{\epsilon}{\cal L}_0+\frac{1}{\sqrt \epsilon}{\cal L}_1+
{\cal L}_2, \label{eq:gav}\\
{\cal L}_{hom}&=\frac{1}{\epsilon^2}{\cal L}_0+\frac{1}{\epsilon}{\cal L}_1+
{\cal L}_2, \label{eq:ghom}
\end{align}
with the understanding that $f_0 \equiv 0$ and $g_1 \equiv 0$ in the case of ${\cal L}_{av}.$
We let $\Omega$ denote the probability space for the pair of Brownian motions $U,V$.

In $\eqref{eq:av}$ (resp. $\eqref{eq:hom}$) the dynamics
for $y$ with $x$ viewed as frozen has solution $\varphi_{x}^{t/\eps}(y(0))$ (resp.
$\varphi_{x}^{t/\epss}(y(0))$). Of course $x$ is not frozen, but since it evolves much
more slowly than $y$, intuition based on freezing $x$ and considering the process
\eqref{eq:fast} is useful in understanding how averaging and homogenization arise for
equations \eqref{eq:av} and \eqref{eq:hom} respectively. Specifically, for \eqref{eq:av}
on timescales long compared with $\epsilon$ and short compared to $1$, $x$ will be
approximately frozen and $y$ will traverse its invariant measure with density
$\rho(y;x)$. We may thus average over this measure and eliminate $y$. Similar ideas hold
for equation \eqref{eq:hom}, but are complicated by the presence of the term
$\eps^{-1}f_0.$ These ideas underly the averaging and homogenization results contained
in the next two subsections.
%
%
\subsection{Averaging}
Define $F:\cX \to \bbR^l$ and $K:\cX \to \bbR^{l \times l}$ by
\begin{equation*}
F(x)  := \int_{\cY} f_1(x,y)\rho(y;x)dy
\end{equation*}
and
\begin{equation*}
K(x)K(x)^T := \int_{\cY} \bigl(\alpha_0(x,y)\alpha_0(x,y)^T +\alpha_1(x,y)
\alpha_1(x,y)^T\bigr)\rho(y;x)dy.
\end{equation*}
Note that $K(x)K(x)^T$ is positive semidefinite and hence $K(x)$ is well defined via, for
example, the Cholesky decomposition.
\begin{theorem}
\label{th:av}
Let Assumptions \ref{ass:a1} hold and let $x(0)=X(0)$. Then
$x \Rightarrow X$ in $C([0,T],{\cal X})$
and $X$ solves the SDE
\begin{equation}
\label{eq:avs}
\frac{dX}{dt}=F(X)+K(X)\frac{dW}{dt},
\end{equation}
where $W$ is ca standard $l$-dimensional Brownian motion.
\end{theorem}
We use the notation $\Omega_0$ to denote the probability space for
the Brownian motion $W$.

\begin{proof}

Consider the Poisson equation
$$-\cL_0 \Xi(y;x)=f_1(x,y)-F(x), \quad \int_{\cY} \rho(y;x) \Xi(y;x)dy=0$$
with unique  solution $\Xi(y;x) \in L^2_{\rho}(\cY;x).$ Applying It\^{o}'s formula to
$\Xi$ we obtain
$$\frac{d\Xi}{dt}=\rato \cL_0 \Xi+\frac{1}{\sqrt \eps}\cL_1 \Xi+\cL_2 \Xi+\rath \nabla_y \Xi \beta
\frac{dV}{dt}+\nabla_x \Xi\alpha_0\frac{dU}{dt}+\nabla_y\Xi\alpha_1 \frac{dV}{dt}.$$
From this we obtain
$$\int_0^t \Bigl( f_1(x(s),y(s))-F(x(s)) \Bigr) ds=e_0(t)$$
where
\begin{align*}
e_0(t)=\sqrt \epsilon \int_0^t \left(\cL_1 \Xi ds+\nabla_y\Xi \beta dV\right)+
&\eps \int_0^t \left(\cL_2 \Xi ds +\nabla_x \Xi \alpha_0 dU+\nabla_y\Xi\alpha_1 dV \right)\\
+&\eps \left(\Xi(y(0);x(0))-\Xi(y(t);x(t))\right).
\end{align*}
Thus, by Assumptions \ref{ass:a1} and the Burkholder-Davis-Gundy inequality,
$$e_0 \to 0\,\,{\mbox{in}}\,\,L^p(C([0,T],{\cal X});\Omega).$$

Hence
\begin{equation*}
x(t)= x(0)+\int_0^t F(x(s))ds+M(t)+e_0(t)
\end{equation*}
with
\begin{equation*}
M(t):= \int_0^t \alpha_0(x(s),y(s))dU(s)+\int_0^t \alpha_1(x(s),y(s))dV(s).
\end{equation*}
The quadratic variation process for $M(t)$ is
$$
\langle M \rangle_t = \int_0^t A(x(s),y(s)) \, ds,
$$
where
$$
A(x,y) = \alpha_0 (x,y) \alpha_0(x,y)^T + \alpha_1(x,y) \alpha_1(x,y)^T.
$$
By use of the Poisson equation technique applied above to show that $f_0(x,y)$
can be approximated by $F(x)$ (its average  against the fast $y$ process),
we can show similarly that
$$\int_0^t A(x(s),y(s))ds=\int_0^t K(x(s)) K(x(s))^T ds+ e_1(t)$$
where, as above,
$$e_1 \to 0\,\,{\mbox{in}}\,\,L^p(C([0,T],{\cal X});\Omega).$$

Let
\begin{align*}
B(t)&=x(0)+\int_0^t F(x(s))ds+  e_0(t),\\
q(t)&=\int_0^t K(x(s))K(x(s))^T ds+ e_1(t).
\end{align*}
Then
$$x(t)=B(t)+M(t),$$
where $M(t)$ and $M(t)M(t)^T-q(t)$ are ${\cal F}_t$ martingales, where ${\cal F}_t$ is
the filtration generated by $\sigma((U(s),V(s)), s \le t).$ Let $C_{c}^{\infty}({\cal
X})$ denote the space of compactly supported $C^{\infty}$ functions. The martingale
problem for
$$
{\mathcal A}=\{(f,K:F \cdot \nabla f + \nabla_x\nabla_x f): f \in C_{c}^{\infty}({\cal
X})\}
$$
is well posed and $x(s),y(s)$ and $X(s)$ are continuous. By $L^2$ convergence of the
$e_i$ to $0$ in $C([0,T],{\cal X})$ we deduce convergence to $0$ in probability, in the
same space. Hence by a slight generalization of Theorem 4.1 in Chapter 7 of
\cite{EthKur86} we deduce the desired result.
\end{proof}

\subsection{Homogenization}

In order for the equations \eqref{eq:hom} to produce a sensible limit
as $\eps \to 0$ it is necessary to impose a condition on $f_0.$ Specifically
we assume the following which, roughly, says that $f_0(x,y)$ averages
to zero against the invariant measure of the fast $y$ process, with $x$ fixed.

\begin{assumptions}
\label{ass:a2} The function $f_0$ satisfies the centering condition
$$\int_{\cY} \rho(y;x)f_0(x,y)dy=0.$$
\end{assumptions}

Let $\Phi(y;x) \in L^2_{\rho}(\cY;x)$ be the solution of the equation
\begin{equation}\label{e:poisson_main}
-\cL_0 \Phi(y;x)=f_0(x,y), \quad \int_{\cY} \rho(y;x) \Phi(y;x)dy=0,
\end{equation}
which is unique by Assumptions \ref{ass:a2}.  Define
\begin{align*}
F_0(x)&:=\int_{\cal Y} ({\cal L}_1\Phi)(x,y)\rho(y;x)dy\\
&=\int_{\cY}\Bigl( \bigl(\nabla_x \Phi f_0\bigr)(x,y)+
\bigl(\nabla_y \Phi g_1 \bigr)(x,y)+ \bigl(\alpha_1 \beta^T:\nabla_y\nabla_x
\Phi\bigr)(x,y) \Bigr)\rho(y;x)dy,\\
F_1(x)&:=\int_{\cY} f_1(x,y) \rho(y;x) dy\,\quad\mbox{and}\\
F(x)&=F_0(x)+F_1(x).
\end{align*}
Also define
\begin{align*}
A_1(x)A_1(x)^T&:=\int_{\cY}\Bigl( \bigl(\nabla_y \Phi \beta+\alpha_1 \bigr)
\bigl(\nabla_y \Phi
\beta+\alpha_1 \bigr)^T \Bigr)(x,y)\rho(y;x)dy,\\
A_0(x)A_0(x)^T&:=\int_{\cY} \alpha_0(x,y) \alpha_0(x,y)^T  \rho(y;x) dy\,\quad \mbox{and}\\
K(x)K(x)^T&=A_0(x)A_0(x)^T+A_1(x)A_1(x)^T.
\end{align*}
Note that $K(x)K(x)^T$ is positive semidefinite by construction so that $K(x)$ is well
defined by, for example, the Cholesky decomposition.
\begin{theorem}
\label{th:hom}
Let Assumptions \ref{ass:a1}, \ref{ass:a2} hold. Then
$x \Rightarrow X$ in $C([0,T],{\cal X})$ and $X$ solves the SDE
\begin{equation}
\label{eq:homs}
\frac{dX}{dt}=F(X)+A(X)\frac{dW}{dt}
\end{equation}
where $W$ is a standard $l$-dimensional Brownian motion.
\end{theorem}

\begin{proof}
We consider three Poisson equations: that for $\Phi$ given above
and
\begin{subequations}
\begin{eqnarray}
-\cL_0 \chi(y;\xi)&=& f_1(x,y)-F_1(x), \quad \int_{\cY} \rho(y;x) \chi(y;x)dy=0,
\label{e:chi}\\
-\cL_0 \Psi(y;\xi)&=& (\cL_1 \Phi)(x,y)-F_0(x), \quad \int_{\cY} \rho(y;x) \Psi(y;x)dy=0.
\label{e:psi}
\end{eqnarray}
\end{subequations}
All of these equations have a unique solution since the right hand sides average to
zero against the density $\rho(y;x)$ by assumption ($\Phi$) or by
construction ($\chi$, $\Psi$).

By the It\^{o} formula we obtain
$$\frac{d\Phi}{dt}=\ratt \cL_0 \Phi+\rato \cL_1 \Phi++\cL_2 \Phi+
\rato \nabla_y \Phi \beta \frac{dV}{dt}+\nabla_x \Phi \alpha_0 \frac{dU}{dt}+
\nabla_x \Phi \alpha_1\frac{dV}{dt}.$$
From this we obtain, using arguments similar to those in the proof of Theorem \ref{th:av},
$$
\rato \int_0^t f_0(x,y)ds=\int_0^t (\cL_1\Phi)(x(s),y(s))ds
+\int_0^t (\nabla_y\Phi \beta)(x(s),y(s))dV(s)+ e_0(t)$$
where
$$e_0(t) \to 0 \;\; {\mbox{in}}\;\; L^p(C([0,T],{\cal X});\Omega)$$
and where, recall, $\Omega$ is the probability space for $(U,V).$ Applying It\^{o}'s
formula to $\chi$, the solution of~\eqref{e:chi}, we may show that
$$\int_0^t \Bigl( f_1(x(s),y(s))-F_1(x(s))\Bigr)ds=e_1(t)$$
where
$$e_1(t) \to 0 \;\; {\mbox{in}} \;\; L^p(C([0,T],\bbR^d);\Omega).$$
Thus
\begin{align*}
x(t)=x(0)&+\int_0^t \bigl(\cL_1\Phi\bigr)(x(s),y(s))ds+\int_0^t F_1(x(s))ds+
\int_0^t \bigl(\nabla_y\Phi \beta\bigr)(x(s),y(s))dV(s)\\
&+\int_0^t \alpha_0(x(s),y(s))dU(s)+\int_0^t \alpha_1(x(s),y(s))dV(s)+e_2(t)
\end{align*}
and
$$e_2(t) \to 0 \;\;  {\mbox{in}} \; \; L^p(C([0,T],{\cal X});\Omega).$$
By applying It\^{o}'s formula to $\Psi$, the solution of~\eqref{e:psi} we obtain
$$\frac{d\Psi}{dt}=\ratt \cL_0 \Psi+\rato \cL_1 \Psi++\cL_2 \Psi+
\rato \nabla_y \Psi \beta \frac{dV}{dt}+\nabla_x \Psi \alpha_0 \frac{dU}{dt}+
\nabla_x \Psi \alpha_1\frac{dV}{dt}$$
From this we obtain
$$\int_0^t \Bigl(\cL_1 \Phi-F_0 \Bigr)(x,y)ds=e_3(t)$$
where
$$e_3(t) \to 0 \;\; {\mbox{in}}  \;\;  L^p(C([0,T],{\cal X});\Omega).$$
Thus
\begin{align*}
x(t)=&x(0)+\int_0^t F(x(s))ds+M(t)+e_4(t)\quad \mbox{and}\\
M(t):=&\int_0^t \alpha_0(x(s),y(s)) \, dU(s)+\bigl(\nabla_y\Phi\beta+\alpha_1\bigr)
(x(s),y(s)) \, dV(s).
\end{align*}
Here
$$e_4 \to 0 \;\; {\mbox{in}}  \;\;  L^p(C([0,T],{\cal X});\Omega).$$
Define
$$A_2(x,y)=\bigl(\nabla_y \Phi \beta+\alpha_1 \bigr) \bigl(\nabla_y \Phi
\beta+\alpha_1 \bigr)^T (x,y)
+\alpha_0(x,y) \alpha_0(x,y)^T .$$
The quadratic variation of $M(t)$ is
$$\langle M \rangle_t = \int_0^t A_2(x(s),y(s)) \, ds.$$
By use of the Poisson equation technique
we can show that
$$\int_0^t A_2(x(s),y(s))ds=\int_0^t K(x(s)) K(x(s))^T ds+ e_5(t)$$
where, as above,
$$e_5 \to 0 \;\; {\mbox{in}}  \;\;  L^p(C([0,T],{\cal X});\Omega).$$
The remainder of the proof proceeds as in Theorem \ref{th:av}.
\end{proof}

\section{Parameter Estimation}\label{sec:par_est}

Recall that $\Omega_0$ is the probability space for $W$.
Imagine that we try to fit data $\{x(t)\}_{t \in [0,T]}$ from \eqref{eq:av}
or \eqref{eq:hom} to a homogenized or averaged
equation of the from \eqref{eq:avs} or \eqref{eq:homs}, but with unknown parameter
$\theta \in \Theta$, where $\Theta$ is an open subset of $\R^k$, in the drift:
\begin{equation}
\label{eq:fitthis}
\frac{dX}{dt}=F(X;\theta)+K(X)\frac{dW}{dt}.
\end{equation}
Suppose that the actual drift compatible with the data is given by $F(X)=F(X;\theta_0).$
We ask whether it is possible to  correctly identify $\theta=\theta_0$ by finding the
{\em maximum likelihood estimator} (MLE) when using a statistical model of the form
\eqref{eq:fitthis}, but given data from \eqref{eq:av} or \eqref{eq:hom}. Recall that the
averaging and homogenization techniques from the previous section show that $x(t)$ from
\eqref{eq:av} and \eqref{eq:hom} converges weakly to the solution of an equation of the
form \eqref{eq:fitthis}. We make the following assumptions concerning the model equations
\eqref{eq:fitthis} which will be used to fit the data.
\begin{assumptions} \label{ass:a0}
We assume that $K$ is uniformly positive-definite on ${\cal X}.$
We also assume that \eqref{eq:fitthis} is ergodic with invariant measure
$\nu(dx)=\pi(x)dx$ at $\theta=\theta_0$ and that
\begin{equation}\label{e:a_infty}
A_{\infty}:=\int_{\cal X} \left(K(x)^{-1} F(x) \otimes K(x)^{-1} F(x)\right) \pi(x)dx
\end{equation}
is invertible.
\end{assumptions}

Given data $\{z(t)\}_{t \in [0,T]}$, the log likelihood function for $\theta$ satisfying
\eqref{eq:fitthis} is given by
\begin{equation}
\label{eq:like} \bL(\theta;z)=\int_0^T \langle F(z;\theta), dz
\rangle_{a(z)}-\frac12\int_0^T |F(z;\theta)|_{a(z)}^2 dt,
\end{equation}
where
$$\langle p, q \rangle_{a(z)}= \langle K(z)^{-1} p, K(z)^{-1} q \rangle.$$
To be precise
$$\frac{d\bbP}{d\bbP_0}=\exp \left({\bL}(\theta;X)\right)$$
where $\bbP$ is the path space measure for \eqref{eq:fitthis} and $\bbP_0$ the pathspace
measure for \eqref{eq:fitthis} with $F \equiv 0$~\cite{Rao99}. The MLE is
\begin{equation}\label{e:argmax}
\htheta=\mbox{argmax}_{\theta} \bL(\theta;z).
\end{equation}
As a preliminary to understanding the effect of using multiscale
data,
we start by exhibiting an underlying property of the log-likelihood
when confronted with data from the model \eqref{eq:fitthis} itself. The following theorem
shows that, in this case: (i) in the limit $T \to \infty$ the
log-likelihood is asymptotically independent of the particular sample
path of \eqref{eq:fitthis} chosen -- it depends only on the invariant measure $\pi$; (ii)
as a consequence we see that, asymptotically, time-ordering of the data is irrelevant to
parameter estimation; (iii) under some additional assumptions, the
large $T$ expression also shows that choosing data from the model \eqref{eq:fitthis}
leads to the correct estimation of drift parameters, in the limit $T\to \infty.$

\begin{theorem} \label{th:basic}
Let Assumptions \ref{ass:a0} hold and let
$\{X(t)\}_{t \in [0,T]}$ be a
sample path of \eqref{eq:fitthis} with $\theta=\theta_0.$ Then,
in $L^2(\Omega_0)$ and almost surely with respect to $X(0)$,
$$\lim_{T \to \infty}\frac{2}{T}\bL(\theta;X)
=\int_{\cX}|F(X;\theta_0)|_{a(X)}^2 \pi(X) dX-\int_{\cX} |F(X;\theta)-
F(X;\theta_0)|_{a(X)}^2 \pi(X)dX.$$
This expression is maximized by choosing $\htheta=\theta_0$, in the limit
$T \to \infty.$
\end{theorem}

\begin{proof}
By Lemmas \ref{l:A1} and \ref{l:A2} in the appendix we deduce that, with all limits
in $L^2(\Omega)$,
\begin{align*}
\lim_{T \to \infty}\frac{1}{T}\bL(\theta;X)
&=\lim_{T \to \infty}\Bigl(\frac{1}{T}\int_{0}^T \langle F(X;\theta), F(X;\theta_0) \rangle_{a(X)}dt\\
&\quad+ \frac{1}{T}\int_0^T \langle F(X;\theta), K(X) dW \rangle_{a(X)} dt
-\frac{1}{2T}\int_{0}^T |F(X;\theta)|_{a(X)}^2 dt\Bigr)\\
&=\int_{\cX} \langle F(X;\theta), F(X;\theta_0) \rangle_{a(X)}\pi(X)dX-
\frac12\int_{\cX} |F(X;\theta)|_{a(X)}^2 \pi(X)dX.
\end{align*}
Completing the square provides the proof.
\end{proof}

In the particular case where the parameter $\theta$ appears linearly in the drift it can
be viewed as an $\bbR^{l \times l}$ matrix $\Theta$ and
\begin{equation}
F(X;\theta)=\Theta F(X) \label{eq:this}
\end{equation}
The correct value for $\Theta$ is thus the $\bbR^{l \times l}$ identity matrix $I$. The
maximum likelihood estimator is
\begin{equation}
\label{eq:Th}
\oT(z;T)=A(z;T)^{-1}B(z;T)
\end{equation}
where
\begin{align*}
A(z;T)&=\frac{1}{T}\int_0^T K(z)^{-1}F(z) \otimes K(z)^{-1}F(z) \, dt,\\
B(z;T)&=\frac{1}{T}\int_0^T K(z)^{-1}dz \otimes K(z)^{-1}F(z);
\end{align*}
if $A(z;T)$ is not invertible then we set $\oT(z;T)=0.$
A result closely related to Theorem \ref{th:basic} is the following
\footnote{The proof is standard and we outline it only for comparison with
the situation in the next subsection
where data from a multiscale model is employed.}:

\begin{theorem} \label{th:basic2} Let Assumptions \ref{ass:a0} hold
and let $\{X(t)\}_{t \in [0,T]}$ be a sample path of \eqref{eq:fitthis}
with $\theta=\theta_0$ so that $F(X;\theta)=F(X)$.
Then
$$\lim_{T \to \infty} \oT(X;T)=I$$
in probability.
\end{theorem}

\begin{proof}
We observe that
\begin{equation*}
B(X;T) = A(X;T)+J_1
\end{equation*}
where
\begin{equation*}
J_1 = \frac{1}{T}\int_0^T dW \otimes K(X)^{-1} F(X)
\end{equation*}
and where $\bbE |J_1|^2={\cal O}(1/T)$ by Lemma \ref{l:A1}.
By ergodicity, and Lemma \ref{l:A2}, we have that
$$A(X;T)=A_{\infty}+J_2$$
where $\bbE |J_2|^2={\cal O}(1/T)$ and $A_{\infty}$ is given by~\eqref{e:a_infty}. By
Assumption~\ref{ass:a0} and for $T$ sufficiently large, $A(z;T)$ is invertible and we
have
$$\oT(X;T)=I+(A_{\infty}+J_2)^{-1}J_1$$
and the result follows.
\end{proof}

\begin{remark}
The invertibility of $A_\infty$ is necessary in order to be able to successfully estimate the drift of the linear system. \end{remark}
In order to prove an analogue of Theorem~\ref{th:basic2} when the drift depends
nonlinearly on the parameter $\theta$ we need to make additional assumptions.
\begin{assumptions}
\label{ass:id}
\begin{itemize}
\item{ We assume that
\begin{equation}
\inf_{|u|>\delta}\int_{\cX} |F(X;\theta_0+u)-F(X;\theta_0)|_{a(X)}^2
\pi(X)dX>\kappa(\delta)>0,\ \forall\delta>0.
\label{eq:id}
\end{equation}}
When~\eqref{eq:id} holds we will say that the system is identifiable.
\item{There exist an $\alpha>0$ and $\hat{F}:\cX\rightarrow \R,$ square integrable with respect to the invariant measure, i.e. $\int_\cX
\hat{F}(X)^2\pi(X)dX<\infty$, such that
\begin{equation}
\label{ass:Holder}
|F(X;\theta)-F(X;\theta^\prime)|_{a(X)}\leq |\theta-\theta^\prime|^\alpha \hat{F}(X)
\end{equation}}
\end{itemize}
\end{assumptions}

Under the above assumption we can prove convergence of the MLE to the correct
value $\theta_0$.
\begin{theorem}
Suppose that Assumptions \ref{ass:a0} and \ref{ass:id} hold. If, in addition, the parameter space $\Theta$ is compact, then
\[ \lim_{T\rightarrow\infty}\hat{\theta}(X;T)=\theta_0\]
in probability.
\end{theorem}
\begin{proof}
It is a straightforward application of the results in \cite{vanZanten01}.
\end{proof}


We now ask whether the likelihood behaves similarly when confronted with
data $\{x(t)\}$ from the underlying multiscale systems \eqref{eq:av} or \eqref{eq:hom}.
To address this issue we make the following natural assumptions regarding the
invariant measure for these underlying multiscale systems.
\begin{assumptions}
\label{ass:a3}
\begin{itemize}
\item The fast/slow SDE \eqref{eq:av} (resp. \eqref{eq:hom}) is ergodic with
invariant measure $\mu^{\eps}(dx dy)$ which is absolutely continuous with
respect to the Lebesgue measure on $\mathcal{X} \times \cal Y$ with smooth
density $\rho^\eps(x,y)$.
\item The limiting SDE~\eqref{eq:avs} or \eqref{eq:homs} is ergodic with
invariant measure $\nu(dx)$ which is absolutely continuous with respect to
the Lebesgue measure on $\cal X$ with smooth density $\pi(x)$.
\item The measure $\mu^{\eps}(dx dy) = \rho^{\eps}(x,y) dx dy$ converges weakly
to the measure $\mu(dx dy) = \pi(x)\rho(y;x)dx dy$ where $\rho(y;x)$ is the invariant
density of the fast process \eqref{eq:fast} given in Assumption \ref{ass:a1} and $\pi(x)$
is the invariant density for \eqref{eq:avs} (resp. \eqref{eq:homs}).
\item The invariant measure $\mu^\eps(dx dy) = \rho^\eps (x,y) dx dy$ satisfies a
Poincar\'{e} inequality with a constant independent of $\eps$: there exists a constant
$C_p$ independent of $\eps$ such that for every mean zero $H^1(\cX \times \cY;
\mu^\eps(dx dy))$ function $f$ we have that
\begin{equation}\label{e:poincare}
\|f \| \leq C_p \|\nabla f \|
\end{equation}
where $\nabla $ represents the gradient with respect to $(x^T, \, y^T)^T$ and $\|\cdot
\|$ denotes the $L^2(\cX \times \cY; \mu^\eps(dx dy))$ norm.
\end{itemize}
\end{assumptions}
We also need to assume that the fast/slow SDEs \eqref{eq:av} and
\eqref{eq:hom} are uniformly elliptic.
\begin{assumption}\label{assump:a4}
Define the matrix field $\Sigma=\gamma\gamma^T$ where
\begin{eqnarray*}
\gamma=\left(
\begin{array}{cc}
\alpha_0 & \alpha_1 \\
0 & \frac{1}{\eps}\beta
\end{array}
\right).
\end{eqnarray*}
Then there is $C_{\gamma}>0$, independent of $\epsilon \to 0$ such that
$$\langle \xi, \Sigma(x,y)\xi \rangle \ge C_{\gamma}|\xi|^2 \quad \forall
(x,y) \in \cX \times \cY, \xi \in \bbR^d.$$
\end{assumption}
\begin{remark}
It is straightforward to show that, when $\mathcal{X} = \T^{\ell}, \;\mathcal{Y} = \T^{\ell - d}$,
Assumptions~\ref{ass:a3} follow from  Assumption~\ref{assump:a4}, using properties of
periodic functions~\cite{PavlSt06}, together with the compactness of the state space.
When $\mathcal{X} = \R^{\ell}, \;\mathcal{Y} = \R^{\ell - d}$ more work is needed in
order to prove that the invariant measure satisfies Poincar\'{e}'s inequality with an
$\eps$ independent constant, since this, essentially, requires to prove that the
generator of the fast/slow system has an $\eps$-independent spectral gap. In this case
where the fast/slow system has a gradient structure with a smooth potential $V(x,y)$,
then simple criteria on the potential have been derived that facilitate determination of
whether or not the invariant measure satisfies the Poincar\'{e} inequality. We refer
to~\cite{Vil04HPI, BakryCattiauxGuillin08} and the references therein for more details.
\end{remark}
%
\subsection{Averaging}

We now ask what happens when the MLE for the averaged equation \eqref{eq:fitthis} is
confronted with data from the original multiscale equation \eqref{eq:av}. The following
result shows that, in this case, the estimator will behave well, for large time and small
$\epsilon$. Large time is always required for convergence of drift parameter estimation,
even when model and data match. In the limit $\epsilon \to 0,$ $X(t)$ from
\eqref{eq:fitthis} approximates $x(t)$ from \eqref{eq:av}.

\begin{theorem}
Let Assumptions \ref{ass:a1}, \ref{ass:a0}, \ref{ass:a3} and~\ref{assump:a4}   hold. Let $\{x(t)\}_{t \in
[0,T]}$ be a sample path of \eqref{eq:av} and $\{X(t)\}_{t \in [0,T]}$ a sample path of
$\eqref{eq:fitthis}$ at $\theta=\theta_0.$ Then the following limits, to be interpreted
in $L^2(\Omega)$ and $L^2(\Omega_0)$ respectively, and almost surely with respect to
$x(0), y(0), X(0),$ are identical: \label{th:avpe}
$$\lim_{\eps \to 0}\lim_{T \to \infty}\frac{1}{T}\bL(\theta;x)=\lim_{T
\to \infty}\frac{1}{T}\bL(\theta;X).$$
\end{theorem}

\begin{proof}
We start by observing that, by Lemma \ref{l:A2} and Assumptions \ref{ass:a3},
\begin{align*}
\lim_{\eps \to 0}\lim_{T \to \infty} \frac{1}{T}\int_0^T
|F(x;\theta)|_{a(x)}^2 dt=&
\lim_{\eps \to 0}\int_{\cX \times \cY}|F(x;\theta)|_{a(x)}^2 \rho^{\eps}(x,y)dxdy\\
=&\int_{\cX \times \cY}|F(x;\theta)|_{a(x)}^2 \pi(x)\rho(y;x)dxdy\\
=&\int_{\cX}|F(x;\theta)|_{a(x)}^2\pi(x)dx,
\end{align*}
where the limits are in $L^2(\Omega).$ Now, from Equation~\eqref{eq:av} it follows that
\begin{align*}
\frac{1}{T}\int_0^T \langle F(x;\theta),dx \rangle_{a(x)}&=\frac{1}{T}
\int_0^T \langle F(x;\theta), f_1(x,y) \rangle_{a(x)}dt\\
&\quad+\frac{1}{T}\int_0^T \langle F(x;\theta),\alpha_0(x,y)dU \rangle_{a(x)}
+\frac{1}{T}\int_0^T \langle F(x;\theta),\alpha_1(x,y)dV \rangle_{a(x)}.
\end{align*}
The last two integrals tend to zero in $L^2(\Omega)$ as $T \to \infty$ by
Lemma \ref{l:A1}.
In order to analyze the first integral on the right hand
side we consider solution of the Poisson equation
$$-\cL_0 \Lambda=\langle F(x;\theta), f_1(x,y)-F(x;\theta_0)\rangle_{a(x)}, \quad \int_{\cY} \rho(y;\xi) \Lambda(y)dy=0.$$
This has a unique  solution $\Lambda(y;x) \in L^2_{\rho}(\cY;x)$
by construction of $F$.

Applying It\^{o}'s formula to $\Lambda$ gives
$$\frac{d\Lambda}{dt}=\rato \cL_0 \Lambda+\frac{1}{\sqrt \epsilon}
\cL_1 \Lambda+\cL_2 \Lambda+\rath \nabla_y \Lambda\beta \frac{dV}{dt}
+\nabla_x \Lambda \alpha_0 \frac{dU}{dt}
+\nabla_x \Lambda \alpha_1 \frac{dV}{dt}$$
which shows that
\begin{align*}
\frac{1}{T}\int_0^T \langle F(x;\theta),f_1(x,y) \rangle_{a(x)}dt=&\frac{1}{T}
\int_0^T \langle F(x;\theta), F(x;\theta_0) \rangle_{a(x)}dt\\
&+\frac{\eps}{T}\int_0^T \bigl(\cL_2 \Lambda\bigr)(x(t),y(t))dt-\frac{\eps}{T}
\Bigl(\Lambda(x(T),y(T))-\Lambda(x(0),y(0))\Bigr)\\
+&\frac{1}{T}\int_0^T \epsr\Bigl(\nabla_y \Lambda \beta)(x(t),y(t))dV(t)+
\left(\cL_1 \Lambda \right)(x(t),y(t)) dt \Bigr)\\
+&\frac{1}{T}\int_0^T \eps\Bigl(\nabla_x \Lambda \alpha_0)(x(t),y(t))dU(t)
+\nabla_y \Lambda \alpha_1)(x(t),y(t))dV(t) \Bigl).
\end{align*}
The stochastic integrals tend to zero in $L^2(\Omega)$ as $T \to \infty$.
By assumption $\Lambda$ is bounded.
Furthermore, in $L^2(\Omega)$,
$$
\frac{1}{T}\int_0^T \bigl(\cL_i \Lambda\bigr)(x(t),y(t))dt \to \int_{\cX \times \cY}
\bigl(\cL_i \Lambda\bigr)(x,y)\rho(y;x)dy, \;\; i=1,2.
$$
Hence we deduce that
\begin{align*}
\lim_{\eps \to 0}\lim_{T \to \infty}
\frac{1}{T}\int_0^T \langle F(x;\theta),f_0(x,y) \rangle_{a(x)}dt=&
\lim_{\eps \to 0}\lim_{T \to \infty}
\frac{1}{T}\int_0^T \langle F(x;\theta), F(x;\theta_0)
\rangle_{a(x)}dt\\
=&\lim_{\eps \to 0}
\int_{\cX \times \cY} \langle F(x;\theta), F(x;\theta_0)
\rangle_{a(x)}\rho^{\eps}(x,y)dxdy\\
=&\int_{\cX} \langle F(x;\theta), F(x;\theta_0) \rangle \pi(x)dx.
\end{align*}
The result follows.
\end{proof}

In the particular case of linear parameter dependence, when the MLE is given
by~\eqref{eq:Th}
we have the following result, showing that the MLE recovers the correct answer from high
frequency data compatible with the statistical model in an appropriate asymptotic limit.

\begin{theorem}
Let Assumptions \ref{ass:a1}, \ref{ass:a0}, \ref{ass:a3} and \ref{ass:Holder} hold.
Assume that $F(X;\theta)$ is given by \eqref{eq:this}. Let $\{x(t)\}_{t \in [0,T]}$ be a
sample path of \eqref{eq:av}. Then ${\widehat \theta}$ given by \eqref{eq:Th} satisfies
$$\lim_{\epsilon \to 0}\lim_{T \to \infty} \oT(x;T)=I$$
in probability.
\end{theorem}

\begin{proof}
Using equation \eqref{eq:av} we find that
\begin{align*}
B(x;T)&=A(x;T)+J_3+J_4, \quad \mbox{where}\\
J_3&=\frac{1}{T}\int_0^T K(x)^{-1}\left(f_1(x,y)-F(x)\right) \otimes K(x)^{-1} F(x)dt,\\
J_4&=\frac{1}{T}\int_0^T K(x)^{-1}\left(\alpha_0(x,y)dU+\alpha_1(x,y)dV\right) \otimes K(x)^{-1} F(x).
\end{align*}
Here, for fixed $\epsilon>0$,
$\bbE |J_4|^2={\cal O}(1/T)$ by Lemma \ref{l:A1} and
$$\lim_{\epsilon \to 0}\lim_{T \to \infty}\bbE |J_3|^2=0$$
by use of the Poisson equation technique.
By ergodicity, and Lemma \ref{l:A2}, we have that
$$A(x;T)=A_{\infty,\epsilon}+J_5$$
where
\begin{equation*}
A_{\infty,\epsilon} :=  \int_{{\cal X} \times {\cal Y}} \left(K(x)^{-1} F(x)
\otimes K(x)^{-1} F(x)\right) \rho^{\epsilon}(x,y)dxdy,
\end{equation*}
with
\begin{equation*}
\lim_{\epsilon \to 0} A_{\infty,\epsilon}=A_{\infty}
\end{equation*}
and, for fixed $\epsilon>0$, $\bbE |J_5|^2={\cal O}(1/T).$

Thus by Assumption~\ref{ass:a0} $A(x;T)$ is invertible for $T$ sufficiently
large, and $\epsilon$ sufficiently small, so that
$$\oT(X;T)=I+(A_{\infty,\epsilon}+J_5)^{-1}\left(J_3+J_4\right).$$
The result follows.
\end{proof}

We would like to show that this also holds for the general case, i.e. if
\[
\hat{\theta}(x;T) := \arg\max_\theta\bL(\theta;x)
\]
then
\[\lim_{\epsilon \to 0}\lim_{T \to \infty} \hat{\theta}(x;T)=\theta_0,\
{\rm in\ probability}.\]
In fact, the following theorem is true for every $\epsilon>0$.
\begin{theorem}\label{thm:mle_av}
Let Assumptions \ref{ass:a1}, \ref{ass:a0}, \ref{ass:id}, \ref{ass:a3} and
\ref{assump:a4} hold and assume that $\theta \in \Theta$, a compact set. Let $\{x(t)\}_{t
\in [0,T]}$ be a sample path of \eqref{eq:av} at $\theta=\theta_0.$ Assume furthermore
that that the marginal of the invariant measure of~\eqref{eq:av} on $\cal X$ $\pi^\eps(x)
dx = \Big( \int_\cY \rho^\epsilon(x,y) dy \Big) dx$ is absolutely continuous with respect
to the invariant measure of the limiting SDE $\pi(x)dx$.  Then, for every $\epsilon>0$,
\begin{equation*}
\lim_{T \to \infty} \hat{\theta}(x;T)=\theta_0,\ {\rm in\ probability}.
\end{equation*}
\end{theorem}
\begin{proof}
Let $g_T(\omega,\theta) := \frac{1}{T}\bL(\theta;x)$ and
\[ g_\infty(\theta) := \int_{\cX \times \cY} \left(\langle F(x;\theta), F(x;\theta_0)
\rangle_{a(x)}-\frac{1}{2}|F(x;\theta)|^2_{a(x)}\right)\rho^{\eps}(x,y)dxdy.
\]
It is straightforward to see that
\[ \arg\max_\theta g_\infty(\theta) = \theta_0 \]
by completing the square. We apply Lemma \ref{consistency}, replacing $\eps$ by
$\frac{1}{T}$, $g_\epsilon$ by $g_T$ and $g_0$ by $g_\infty$. The result follows,
provided that conditions (\ref{ass 1a}), (\ref{ass 1b}) and (\ref{ass 2}) are satisfied.
Condition (\ref{ass 1a}) follows from Theorem \ref{th:avpe}. The identifiability
condition (\ref{ass 2}) follows from Assumptions \ref{ass:id} and the absolute continuity
of $\pi^\eps (x) dx = \Big(\int_\cY \rho^\epsilon(x,y) dy \Big) dx$ with respect to
$\pi(x)dx$. Finally, we can verify that (\ref{ass 1b}) holds, following the proof in
\cite{vanZanten01} and using the fact that functions $f_1$, $\alpha_0$ and $\alpha_1$ are
uniformly bounded.
\end{proof}

%
%
\subsection{Homogenization}
We now ask what happens when the MLE for the homogenized
equation \eqref{eq:fitthis} is confronted with data from the
multiscale equation \eqref{eq:hom}, which homogenizes to
give \eqref{eq:fitthis}. The situation differs substantially from
the case where data is taken from the multiscale equations
\eqref{eq:av} which averages to give \eqref{eq:fitthis}:
the two likelihoods are not identical in the large $T$ limit.

In order to state the main result of this subsection we need to introduce the Poisson equation
\begin{equation}
\label{eq:peg}
-\cL_0 \Gamma=\langle F(x;\theta),f_0(x,y)\rangle_{a(x)}, \quad \int_{\cY} \rho(y;\xi) \Gamma(y;x)dy=0
\end{equation}
which has a unique  solution $\Gamma(y;x) \in L^2_{\rho}(\cY;x).$
Note that
$$\Gamma=\langle F(x;\theta),\Phi(x,y) \rangle_{a(x)},$$
where $\Phi$ solves \eqref{e:poisson_main}. Define
\begin{equation}\label{e:Einfty}
E_{\infty}(\theta)=\int_{\cX \times \cY} \Bigl(\cL_1 \Gamma(x,y)- \langle F(x;\theta),
\bigl(\cL_1 \Phi(x,y)\bigr) \rangle_{a(x)}\Bigr)\pi(x)\rho(y;x)dx dy.
\end{equation}
The following  theorem shows that the correct limit of the log likelihood is
not obtained unless $E_{\infty}=0$, something which will not be
true in general. However in the case where $f_0, g_1 \equiv 0$ we do obtain
$E_{\infty}=0$ and in this case we recover the averaging situation
covered in the Theorems \ref{th:av} and Theorem \ref{th:avpe}
(with $\eps$ replaced by $\eps^2$).

\begin{theorem}
\label{th:hompe} Let Assumptions \ref{ass:a1}, \ref{ass:a2}, \ref{ass:a0}, \ref{ass:a3}
and \ref{ass:Holder} hold. Let $\{x(t)\}_{t \in [0,T]}$ be a sample path of
\eqref{eq:hom} and $\{X(t)\}_{t \in [0,T]}$ a sample path of $\eqref{eq:fitthis}$ at
$\theta=\theta_0.$ Then the following limits, to be interpreted in $L^2(\Omega)$ and
$L^2(\Omega_0)$ respectively, and almost surely with respect to $x(0), y(0), X(0)$, are
identical:
$$\lim_{\eps \to 0}\lim_{T \to \infty}\frac{1}{T}\bL(\theta;x)=\lim_{T
\to \infty}\frac{1}{T}\bL(\theta;X)+E_{\infty}(\theta).$$
\end{theorem}

\begin{proof} As in the averaging case of Theorem \ref{th:avpe} we have
$$\lim_{\eps \to 0}\lim_{T \to \infty} \frac{1}{T}\int_0^T
|F(x;\theta)|_{a(x)}^2 dt=\int_{\cX}|F(x;\theta)|_{a(x)}^2\pi(x)dx.$$
Now
$$\frac{1}{T}\int_0^T \langle F(x;\theta), dx \rangle_{a(x)}=I_1+I_2+I_3$$
where
\begin{align*}
I_1&=\frac{1}{\eps T}\int_0^T \langle F(x;\theta), f_0(x,y) \rangle_{a(x)}dt,\\
I_2&=\frac{1}{T}\int_0^T \langle F(x;\theta), f_1(x,y) \rangle_{a(x)}dt,\\
I_3&=\frac{1}{T}\int_0^T \langle F(x;\theta), \alpha_0(x,y)dU+\alpha_1(x,y)dV \rangle_{a(x)}.
\end{align*}
Now $I_3$ is ${\cal O}(1/\sqrt{T})$ in $L^2(\Omega)$ by Lemma \ref{l:A1}.
Techniques similar to
those used in the proof of Theorem \ref{th:avpe} show that
$$\lim_{\eps \to 0} \lim_{T \to \infty}I_2 \to \int_{\cX} \langle F(x;\theta),
F_1(x;\theta_0) \rangle_{a(x)} \pi(dx).$$ Now consider $I_1$. Applying It\^{o}'s formula to the
solution $\Gamma$ of the Poisson equation \eqref{eq:peg}, we obtain
$$\frac{d\Gamma}{dt}=\ratt \cL_0 \Gamma+\rato \cL_1 \Gamma+\cL_2 \Gamma+
\rato \nabla_y \Gamma \beta \frac{dV}{dt}+
\nabla_x \Gamma \alpha_0 \frac{dU}{dt}+ \nabla_x \Gamma \alpha_1\frac{dV}{dt}.$$
From this we deduce that
$$\frac{1}{\eps T}\int_0^T \langle F(x;\theta), f_0(x,y) \rangle dt=
\frac{1}{T}\int_0^T \Bigl(\cL_1 \Gamma\Bigr) dt+I_4$$
where
$$\lim_{\eps \to 0}\lim_{T \to \infty}I_4=0.$$
Thus
$$I_1=\frac{1}{\eps T}\int_0^T \langle F(x;\theta), f_0(x,y) \rangle dt=I_4+I_5+I_6$$
where, in $L^2(\Omega)$,
\begin{align*}
I_5&=\frac{1}{T}\int_0^T \langle F(x;\theta), \bigl(\cL_1 \Phi(x,y)\bigr) \rangle_{a(x)}dt,\\
I_6&=\frac{1}{T}\int_0^T \Bigl(\cL_1 \Gamma(x,y)- \langle F(x;\theta),
\bigl(\cL_1 \Phi(x,y)\bigr) \rangle_{a(x)}\Bigr)dt.
\end{align*}

By the methods used in the proof of Theorem \ref{th:avpe}
we deduce that
$$\lim_{\eps \to 0} \lim_{T \to \infty} I_5 \to \int_{\cX} \langle F(x;\theta),
F_0(x;\theta_0) \rangle_{a(x)} \pi(x)dx.$$
Putting together all the estimates we deduce that, in $L^2$,
\begin{align*}
\lim_{\epsilon \to 0}\lim_{T \to \infty} \frac{1}{T}\bL(x;\theta)
&= \lim_{T \to \infty}\bL(X;\theta)+\lim_{\eps \to 0}\lim_{T \to \infty} I_6\\
&= \lim_{T \to \infty}\bL(X;\theta)+E_{\infty}(\theta).
\end{align*}
\end{proof}
\section{Subsampling}\label{sec:subsamp}

In the previous section we studied the behavior of estimators when confronted with
multiscale data. The data is such that, in an appropriate asymptotic limit $\eps \to 0$,
it behaves weakly as if it comes from a single scale equation in the form of the
statistical model. By considering the behavior of continuous time estimators in the limit
of large time, followed by taking $\eps \to 0$, we studied the behavior of estimators
which do not subsample the data. We showed that in the averaging set-up this did not
cause a problem -- the likelihood behaves as if confronted with data from the statistical
model itself; but in the homogenization set-up the likelihood function was asymptotically
biased for large time. In this section we show that subsampling the data can overcome
this issue, provided the subsampling rate is chosen appropriately.

In the following we use $\bbE^{\pi}$ to denote expectation on $\cX$
with respect to measure with density $\pi$ and
$\bbE^{\rho^\eps}$ to denote expectation on $\cX \times \cY$
with respect to measure with density $\rho^\eps.$ Recall
that, by Assumption~\ref{ass:a3} the latter measure has weak limit with density $\pi(x)\rho(y;x).$
Let $\Omega'=\Omega \times {\cal X} \times {\cal Y}$ and consider
the probability measure induced on paths $x,y$ solving \eqref{eq:hom}
by choosing initial conditions distributed according to the  measure
$\pi(x)\rho(y;x)dxdy.$ With expectation $\bbE$ under this measure
we will also use the notation
$$
\| \cdot \|_p := \left( \E | \cdot|^p \right)^{1/p}.
$$

We define the discrete log likelihood function found from applying the likelihood
principle to the Euler-Marayama approximation of the statistical model
\eqref{eq:fitthis}. Let $z=\{z_n\}_{n=0}^{N-1}$ denote a time series in $\cX$. We obtain
the likelihood
$$
\bL^{\delta,N}(\theta;z) = \sum_{n=0}^{N-1} \langle F(z_n;\theta), z_{n+1} - z_n
\rangle_{a(z_n)}  - \frac{1}{2} \sum_{n=0}^{N-1} |F(z_n ; \theta)|_{a(z_n)}^2 \delta.
$$
Let $x_n=x(n\delta)$, noting that $x(t)$ depends on $\epsilon$, and set
$x=\{x_n\}_{n=0}^{N-1}$. The basic theorem in this section proves convergence of the log
likelihood function, provided that we subsample (i.e. choose $\delta$) at an appropriate
$\epsilon$-dependent rate. We state and prove the theorem, relying on a pair of
intuitively reasonable propositions which we then prove at the end of the section.
\begin{theorem}
\label{th:main1} Let Assumptions \ref{ass:a1}, \ref{ass:a2}, \ref{ass:a0}, \ref{ass:a3}
and \ref{ass:Holder} hold. Let $\{x(t)\}_{t \in [0,T]}$ be a sample path of
\eqref{eq:hom} and $X(t)$ a sample path of $\eqref{eq:fitthis}$ at $\theta=\theta_0.$ Let
$\delta = \eps^{\alpha}$ with $\alpha \in (0,1)$ and let $N = [\eps^{-\gamma} ]$ with
$\gamma > \alpha$. Then the following limits, to be interpreted in $L^2(\Omega')$ and
$L^2(\Omega_0)$ respectively, and almost surely with respect to $X(0)$, are identical:
\begin{equation}\label{e:logl_subsam}
\lim_{\eps \rightarrow 0}\frac{1}{N \delta} \bL^{N, \delta}(\theta;x)
=\lim_{T \to \infty}\frac{1}{T}\bL(\theta;X).
\end{equation}
\end{theorem}
The proof of this theorem is based on the following two technical results,
whose proofs are presented in the appendix.
\begin{prop}\label{cor:x_estim}
Let $( x(t), \, y(t) )$ be the solution of~\eqref{eq:hom} and assume that
Assumptions~\ref{ass:a1} and~\ref{ass:a2} hold. Then, for $\eps, \, \delta$ sufficiently
small, the increment of the process $x(t)$ can be written in the form
\begin{equation}\label{e:x_estim}
x_{n+1} - x_n = F(x_n;\theta_0) \, \delta + M_n + R(\eps, \delta),
\end{equation}
where $M_n$ denotes the martingale term
$$
M_n  =  \int_{n \delta}^{(n+1) \delta} \left(\nabla_y \Phi\beta + \alpha_0 \right) (x(s),
y(s)) \, d V + \int_{n \delta}^{(n+1) \delta} \alpha_1 (x (s) , y(s)) \, d U
$$
with $\|M_n \|_p \leq C \sqrt{\delta}$ and
$$
\| R(\eps, \delta) \|_p \leq C (\delta^{3/2} + \eps \delta^{\frac12} + \eps).
$$
\end{prop}

\smallskip

\begin{prop} \label{lem:ergodic}
Let $g \in C^1(\cX)$ and let Assumptions \ref{ass:a3} hold. Assume that $\eps$ and $N$
are related as in Theorem \ref{th:main1}. Then
\begin{equation}
\lim_{\eps \to 0} \frac{1}{N}\sum_{n=0}^{N-1} g(x_n)=\bbE^{\pi} g, \label{e:ergodic}
\end{equation}
where the convergence is in $L^2$ with respect to the measure on initial conditions with
density $\pi(x)\rho(y;x).$
\end{prop}

\smallskip
\noindent {\bf Proof of Theorem~\ref{th:main1}.} We define
$$
I_1 (x, \theta) = \sum_{n=0}^{N-1} \langle F(x_n;\theta), x_{n+1} - x_n \rangle_{a(x_n)}
$$
and
$$
I_2 (x) = \frac{1}{2} \sum_{n=0}^{N-1} |F(x_n ; \theta)|_{a(x_n)}^2 \delta.
$$
By Proposition~\ref{lem:ergodic} we have that
$$
\frac{1}{N\delta} I_2(x) \to \frac12 \int_{\cX} |F(x;\theta)|_{a(x)}^2\pi(dx).
$$
We use Proposition~\ref{cor:x_estim} to deduce that
\begin{eqnarray*}
\frac{1}{N \delta} I_1 (x ; \theta) & = & \frac{1}{N \delta} \sum_{n=0}^{N-1}\langle
F(x_n;\theta), F(x_n;\theta_0) \delta + M_n + R(\eps, \delta) \rangle_{a(x_n)}
\\ & = & \frac{1}{N }
\sum_{n=0}^{N-1} \langle F(x_n;\theta), F(x_n;\theta_0) \rangle_{a(x_n)} + \frac{1}{N
\delta} \sum_{n=0}^{N-1} \langle F(x_n), M_n \rangle_{a(x_n)} \\ &&+ \frac{1}{N
\delta} \sum_{n=0}^{N-1} \langle F(x_n), R(\eps, \delta) \rangle_{a(x_n)} \\
&=: & J_1 + J_2 + J_3.
\end{eqnarray*}
Again using Proposition~\ref{lem:ergodic} we have that
$$
J_1 \to \int_{\cX} \langle F(x;\theta), F(x; \theta_0) \rangle_{a(x)} \, \pi(dx).
$$
Furthermore, using the fact that $M_n$ is independent of $x_n$ and has
quadratic variation of order $\delta$ it follows that
\begin{eqnarray*}
\|J_2 \|_2^2 & \leq & \frac{1}{N^2 \delta^2} \sum_{n=0}^{N-1} \E \big| \langle F(x_n ;
\theta), M_n \rangle_{a(x_n)} \big|^2 \\
& \leq &\frac{C}{N \delta}.
\end{eqnarray*}
Here $Q$ is defined to obtain the correct quadratic variation of the $M_n$.
Consequently, and since $\gamma > \alpha$,
$$
\|J_2 \|_2 \leq o(1)
$$
as $\eps \to 0$. Similarly, using martingale moment inequalities~\cite[Eq. (3.25) p.
163]{KSh91} we obtain
$$
\|J_2 \|_p \leq o(1).
$$
Finally, again using Proposition~\ref{cor:x_estim}, we have, for
$q^{-1}+p^{-1}=1$,
\begin{eqnarray*}
\|J_3 \|_p & \leq & \frac{1}{N \delta} \sum_{n=0}^{N-1} \|F(x_n) \|_{q} \| R(\eps,
\delta)\|_{p} \leq C \frac{1}{N \delta} N \Big(\delta^{3/2} + \eps + \eps \delta^{1/2}
\Big) \\ & \leq & o(1),
\end{eqnarray*}
as $\epsilon \to 0$, since we have assumed that $\alpha \in (0,1)$.

We thus have
$$\lim_{\eps \to 0} \frac{1}{N\delta}{\bL}^{N,\delta}(\theta;x)=
\int_{\cX} \langle F(x;\theta), F(x; \theta_0) \rangle_{a(x)} \, \pi(x)dx
-\frac12 \int_{\cal X} |F(x;\theta)|_{a(x)}^2\pi(x)dx.
$$
By completing the square we obtain~\eqref{e:logl_subsam}. \qed

\smallskip

As before, we would like to use this theorem in order to prove the consistency of our estimator. The theory developed in~\cite{vanZanten01} no longer applies because it is based on the assumption that the function we are maximizing
(i.e. the log likelihood function) is a continuous semimartingale, which is not true
for the discrete semimartingale $\cL^{N, \delta}(\theta;x)$. The most difficult part
in proving consistency is to prove that the martingale converges uniformly to zero
(Assumption \ref{ass 1b} in Lemma \ref{consistency}). To avoid this difficulty, we make
some extra assumptions that allow us to get rid of the martingale part:
\begin{assumptions}
\label{ass:h}
\begin{enumerate}
\item{
There exists a function $V:\cX\times\Theta\rightarrow \R$ such that for each
$\theta\in\Theta$,  $V(\cdot,\theta)\in{C}^3(\cX)$ and
\begin{equation}
\label{eq:ass:h}
\nabla V(z;\theta) = \left( K(z)K(z)^T\right)^{-1}F(z;\theta),\ \
\forall z\in\cX,\theta\in\Theta.
\end{equation}}
\item{Define $G:\cX\times\Theta\rightarrow\R$ as follows:
\begin{equation*}
G(z;\theta) := D^2 V(z;\theta) : (K(z)K(z)^T),
\end{equation*}
where $D^2 V$ denotes the Hessian matrix of $V$. Then there exist an $\beta>0$ and
$\hat{G}:\cX\rightarrow \R$ that is square integrable with respect to the invariant
measure, such that
\[|G(z;\theta)-G(z;\theta^\prime)|\leq |\theta-\theta^\prime|^\beta \hat{G}(z).\]}
\end{enumerate}
\end{assumptions}
Suppose that the above assumption is true and $\{X(t)\}_{t \in [0,T]}$ is a sample path of $\eqref{eq:fitthis}$. Then, if we apply It\^o's formula to function $V$, we get that for every $\theta\in\Theta$:
\[
dV(X(t);\theta) =  \langle \nabla V(X(t);\theta),dX(t) \rangle + \frac{1}{2}
G(X(t);\theta)dt.\] But from (\ref{eq:ass:h}) we have that
\begin{eqnarray*}
 \langle \nabla V(X(t);\theta),dX(t) \rangle &=&  \langle
\left( K(X(t))K(X(t))^T\right)^{-1}F(X(t);\theta),dX(t) \rangle = \\
&=& \langle F(X(t);\theta),dX(t) \rangle_{a(X(t))}
\end{eqnarray*}
and thus
\begin{equation*}
\langle F(X(t);\theta),dX(t) \rangle_{a(X(t))} = dV(X(t)) - \frac{1}{2} G(X(t);\theta)dt.
\end{equation*}
Using this identity, we can write the log-likelihood function~\eqref{eq:like} in the form
\[
\bL(\theta;X(t))= \left( V(X(T);\theta)-V(X(0);\theta) \right) -\frac12\int_0^T
\left(|F(X(t);\theta)|_{a(X(t))}^2 + G(X(t);\theta) \right) \, dt.
\]
Using this version of the log-likelihood function , we define
\begin{equation}
\label{eq:newL} \tilde\bL^{N,\delta}(\theta;z) = - \frac{1}{2} \sum_{n=0}^{N-1} \left(
|F(z_n ; \theta)|_{a(z_n)}^2 + G(z_n ; \theta)\right)\delta.
\end{equation}
Now we can prove asymptotic consistency of the MLE, provided that we subsample
at the appropriate sampling rate.
\begin{theorem}\label{thm:subsam}
Let Assumptions \ref{ass:a1}, \ref{ass:a2}, \ref{ass:a0}, \ref{ass:id}, \ref{ass:a3},
\ref{assump:a4} and
\ref{ass:h} hold and assume that $\theta \in \Theta,$ a compact set.
Let $\{x(t)\}_{t \in [0,T]}$ be a sample path of \eqref{eq:hom} at $\theta=\theta_0$. Define
\[ \hat{\theta}(x;\eps) := \arg\max_\theta \tilde\bL^{N,\delta}(\theta;x)\]
with $N$ and $\delta$ as in Theorem \ref{th:main1} above and
$\tilde\bL^{N,\delta}(\theta;x)$ defined in~\eqref{eq:newL}. Then,
\begin{equation*}
\lim_{\eps \to 0} \hat{\theta}(x;\epsilon)=\theta_0,\ {\rm in\ probability}.
\end{equation*}
\end{theorem}
\begin{proof}
We apply Lemma \ref{consistency} with $g_\epsilon(x,\theta) =
\frac{1}{N \delta} \tilde{\cL}^{N, \delta}(\theta;x)$ and $g_0(\theta)$ its limit. Note that
\begin{equation*}
\lim_{\eps \rightarrow 0}\frac{1}{N \delta} \tilde{\bL}^{N, \delta}(\theta;x) =\lim_{T
\to \infty}\frac{1}{T}\bL(\theta;X)
\end{equation*}
by Proposition \ref{lem:ergodic} and the fact that
\[
\lim_{T\rightarrow\infty}\frac1T \left(V(X(T);\theta)-V(X(0);\theta)\right) = 0,
\]
which follows from the ergodicity of $X$. As in Theorem \ref{th:main1}, the limits are
interpreted in $L^2(\Omega')$ and $L^2(\Omega_0)$ respectively, and almost surely with
respect to $X(0)$. As we have already seen, the maximizer of $g_0(\theta)$ is $\theta_0$.
So, Assumption (\ref{ass 1a}) is satisfied. Also, Assumption \ref{ass:id} is equivalent
to (\ref{ass 2}). To prove consistency, we need to prove (\ref{ass 1b}), which can be
viewed as uniform ergodicity. The proof is again similar to that in \cite{vanZanten01}.
First, we note that by Assumptions \ref{ass:id} and \ref{ass:h}, both
$g_\eps(\cdot,\theta)$ and $g_0(\theta)$ are continuous with respect to $\theta$, so it
is sufficient to prove (\ref{ass 1b}) on a countable dense subset $\Theta^\star$ of
$\Theta$. Then, uniform ergodicity follows from \cite[Thm. 6.1.5]{dudley84} , provided
that
\[ N_{[\ ]}\left(\eps,{\mathcal F},\|\cdot\|_{L^1(\pi)}\right)<\infty,\]
i.e. the number of balls of radius $\eps$ with respect to $\|\cdot\|_{L^1(\pi)}$ needed to
cover
\[{\mathcal F}:= \{ |F(z;\theta|^2_{a(z)} + G(z;\theta) ;\ \theta\in\Theta^\star\}\]
is finite. As demonstrated in \cite{vanZanten01}, this follows from the H\"{o}lder
continuity of $|F(z;\theta)|^2_{a(z)}$ and $G(z;\theta)$.
\end{proof}

\section{Examples}\label{sec:examp}

Numerical experiments, illustrating the phenomena studied in this paper, can be found in
the paper \cite{PavlSt06}. The experiments therein are concerned with a particular case
of the general homogenization framework considered in this paper and illustrate the
failure of the MLE when the data is sampled too frequently, and the role of subsampling
to ameliorate this problem. In this section we construct two examples which identify the
term $E_{\infty}$ responsible for the failure of the MLE.

\subsection{Langevin Equation in the High Friction Limit}
We consider the Langevin equation in the high friction limit:\footnote{ We have rescaled
the equation in such a way that we actually consider the small mass, rather than the high
friction limit. In the case where the mass and the friction are scalar quantities the two
scaling limits are equivalent.}
\begin{equation}\label{e:lang}
\eps^2 \frac{d^2 q}{dt^2} = - \nabla_q V(q; \theta) - \frac{dq}{dt} + \sqrt{2
\beta^{-1}}\frac{dW}{dt},
\end{equation}
where $V(q;\theta)$ is a smooth confining potential depending on a parameter
$\theta \in \Theta \subset \R^{\ell}$,\footnote{A standard example is that
of a quadratic potential $V(q;\theta) = \frac{1}{2} q \theta q^T$ where the
parameters to be estimated from time series are the elements of the stiffness
matrix $\theta$.} $\beta$ stands for the inverse temperature
and $W(t)$ is standard Brownian motion on $\R^d$.
We write this equation as a first order system
\begin{equation}\label{e:syst}
\frac{dq}{dt} = \frac{1}{\eps} p, \quad \frac{dp}{dt} = - \frac{1}{\eps} \nabla_q V(q;
\theta) - \frac{1}{\eps^2} p + \sqrt{\frac{2 \beta^{-1}}{\eps^2}} \frac{dW}{dt}.
\end{equation}
In the notation of the general homogenization set-up we have
$(x,y)=(q,p)$ and
$$f_0=p, \;\; f_1=0, \;\; \alpha_0=0,\;\; \alpha_1=0$$
and
$$g_0=-p,\;\; g_1=-\nabla_q V(q), \;\; \beta \mapsto \sqrt{2\beta^{-1}}I.$$
The fast process is simply an Ornstein-Uhlenbeck process with generator
$$
{\cal L}_0=-p \cdot \nabla_p+\beta^{-1} \Delta_p.
$$
The unique square integrable (with respect to the invariant measure of the
OU process) solution of the Poisson equation~\eqref{e:poisson_main} is $\Phi
= p$. Therefore,
$$F_0=-\nabla_q V(q; \theta),\;\; F_1=0,\;\; A_1=\sqrt{2 \beta^{-1}}I.$$
Hence the homogenized equation is\footnote{In this case we can actually prove strong
convergence of $q(t)$ to $X(t)$~\cite{nelson,PavlSt03}.}
\begin{equation}\label{e:homog_grad}
\frac{dX}{dt} = - \nabla V(X;\theta) + \sqrt{2 \beta^{-1}}\frac{dW}{dt}.
\end{equation}
Consider now the parameter estimation problem for "full dynamics"~\eqref{e:lang}
and the "coarse grained" model~\eqref{e:homog_grad}: We are given data from~\eqref{e:lang} and we want to fit it to
equation~\eqref{e:homog_grad}. Theorem~\ref{th:hompe} implies that for this problem the
maximum likelihood estimator is asymptotically
biased.\footnote{Subsampling, at the
rate given in Theorem~\ref{th:main1}, is necessary for the correct estimation
of the parameters in the drift of the homogenized equation~\eqref{e:homog_grad}.}
In fact, in this case we can compute the term $E_{\infty}$, responsible for
the bias and  given in
equation~\eqref{e:Einfty}. We have the following result.
\begin{prop}
Assume that the potential $V(q;\theta) \in C^{\infty}(\R^d)$ is such that $e^{-\beta V(q;\theta)}
\in L^1(\R^d)$ for every $\beta >0$ and all $\theta \in \Theta$. Then error term $E_{\infty}$, eqn.~\eqref{e:Einfty} for the SDE~\eqref{e:lang} is given
by the formula
\begin{equation}\label{e:error_grad}
E_{\infty}(\theta) = -Z^{-1}_V \frac{\beta}{2} \int_{\R^d} |\nabla_q V(q;\theta)|^2  e^{-
\beta V(q;\theta)} \, dq,
\end{equation}
where $Z_V = \int_{\R^d} e^{-\beta V(q;\theta)} \, dq$. In particular, $E_{\infty} <0$.
\end{prop}
\proof
We have that
$$
\cL_1 =p \cdot \nabla_q - \nabla_q V \cdot \nabla_p.
$$
The invariant measure of the process is $\epsilon$-independent and we write
it is
$$
\rho(q,p; \theta) \, dq dp = Z^{-1} e^{- \beta H(p,q; \theta)} \, dq dp.
$$
Furthermore, since the homogenized diffusion matrix is $\sqrt{2 \beta^{-1}}I$,
$$
\langle \cdot , \cdot \rangle_{a(z)} = \frac{\beta}{2} \langle \cdot , \cdot \rangle,
$$
where $\langle \cdot , \cdot \rangle$ stands for the standard Euclidean inner product.
We readily check that
$$
\frac{2}{\beta}{\cal L}_1 \Gamma=\cL_1 \langle -\nabla_q V, p \rangle
= - p \otimes p : D^2_q V(q;\theta) + |\nabla_q V(q;\theta)|^2
$$
and
$$
\frac{2}{\beta}\langle F,{\cal L}_1\Phi \rangle_a=\langle -
\nabla_q V, \cL_1 p \rangle = |\nabla_q V(q;\theta)|^2.
$$
Thus,
\begin{eqnarray*}
E_{\infty}(\theta) &=& -\frac{\beta}{2} \int_{\R^{2 d}} p \otimes p : D^2_q V(q;\theta)
Z^{-1} e^{- \beta H(p,q;\theta)} \, dq dp \\ & = & - \frac{1}{2} \int_{\R^d} \Delta_q
V(q;\theta) Z^{-1}_V e^{- \beta V(q;\theta)} \, dq = - \frac{\beta}{2} \int_{\R^d}
|\nabla_q V(q;\theta)|^2 Z^{-1}_V e^{- \beta V(q;\theta)} \, dq,
\end{eqnarray*}
which is precisely~\eqref{e:error_grad}. \qed
%
\subsection{Motion in a Multiscale Potential}

Consider the equation~\cite{PavlSt06}
\begin{equation}
\label{eq:pot}
\frac{dx}{dt}=-\nabla V^{\eps}(x)+\sqrt{2  \beta^{-1}  }\frac{dW}{dt}
\end{equation}
where
$$V^{\eps}(x)=V(x)+p(x/\eps),$$
where the fluctuating part of the potential $p(\cdot)$ is taken to be a smooth
$1$-periodic function.

Setting $y=x/\eps$ we obtain
\begin{subequations}\label{e:mult_pot}
\begin{eqnarray}
\frac{dx}{dt}&=-\Bigl( \nabla V(x)+\rato\nabla p(y)\Big)+\sqrt{2 \beta^{-1} }\frac{dW}{dt}\\
\frac{dy}{dt}&=-\rato \Bigl( \nabla V(x)+\rato\nabla p(y)\Big)+\rato
\sqrt{2  \beta^{-1} }\frac{dW}{dt}.
\end{eqnarray}
\end{subequations}
In the notation of the general homogenization set-up we have
$$f_0=g_0=-\nabla_y p(y),\;\; f_1=g_1= -\nabla V(x)$$
and
$$\alpha_0=0,\;\; \alpha_1=\beta=\sqrt{2\beta^{-1}}.$$
The fast process has generator
$${\cal L}_0=-\nabla_y p(y) \cdot \nabla_y+ \beta^{-1} \Delta_y.$$
The invariant density is $\rho(y) = Z_p^{-1}\exp(-\beta p(y))$ with $Z_p = \int_{\T^d}
\exp(- \beta p(y)) \, dy.$
The Poisson equation for $\Phi$ is
$${\cal L}_0\Phi(y)=\nabla_y p(y).$$
Notice that $\Phi$ is a function of $y$ only. The homogenized equation is
\begin{equation}\label{e:homog_pot}
\frac{d X}{d t} = -K \nabla V(X) + \sqrt{2 \beta^{-1} K} \frac{d W}{d t}
\end{equation}
where
$$
K = \int_{\T^d} (I + \nabla_y \Phi(y)) (I + \nabla_y \Phi(y))^T \rho(y) \,
dy.
$$
Suppose now that the potential contains parameters, $V = V (x, \theta), \; \theta \in
\Theta \subset \R^{\ell}$. We want to estimate the parameter $\theta$, given data
from~\eqref{eq:pot} and using the homogenized equation
$$
\frac{d X}{d t} = - K\nabla V(X; \theta) + \sqrt{2\beta^{-1} K} \frac{d W}{d t}.
$$
Theorem~\ref{th:hompe} implies that, for this problem, the maximum likelihood
estimator is asymptotically biased and that subsampling at the appropriate
rate is necessary for the accurate estimation of the parameter $\theta$. As
in the example presented in the previous section, we can calculate explicitly
the error term $E_{\infty}$. For simplicity we will consider the problem
in one dimension.
\begin{prop}
Assume that the potential $V(x;\theta) \in C^{\infty}(\R)$ is such that
$e^{-\beta V(x;\theta)} \in L^1(\R)$ for every $\beta >0$ and all $\theta \in \Theta$.
Then error term $E_{\infty}$, eqn.~\eqref{e:Einfty} for the SDE~\eqref{eq:pot} is given
by the formula
\begin{equation}\label{e:error_pot}
E_{\infty}(\theta) = \Big(-1 + \widehat{Z}_p^{-1} Z_p^{-1} \Big)\frac{\beta Z_V^{-1}}{2}
\int_{\R} |\partial_x V|^2 e^{-\beta V(x; \theta)} \, dx .
\end{equation}
where $Z_V = \int_{\R} e^{- \beta V(q;\theta)} \, dq, \, Z_p = \int_0^1 e^{-\beta p(y)}
\, dy \, \widehat{Z}_p = \int_0^1 e^{\beta p(y)} \, dy $. In particular, $E_{\infty} <0$.
\end{prop}

\proof
Equations~\eqref{e:mult_pot} in one dimension become
\begin{subequations}
\begin{eqnarray}\label{e:x_y}
\dot{x} &=& -\partial_x V(x;\theta) - \frac{1}{\eps} \partial_y p(y) + \sqrt{2
\beta^{-1}}
\dot{W}, \\
\dot{y} & = & - \frac{1}{\eps} \partial_x V(x;\theta) - \frac{1}{\eps^2} \partial_y p(y)
+ \frac{2 \beta^{-1}}{\eps^2} \dot{W}.
\end{eqnarray}
\end{subequations}
The invariant measure of this system is (notice that it is independent of $\eps$)
$$
\rho(y,x;\theta) \, dx dy = Z_V^{-1}(\theta) Z_p^{-1} e^{-\beta V(x;\theta) - \beta p(y)}
\, dx dy.
$$
The homogenized equation is
$$
\dot{X} = -K \partial_x V(x;\theta) + \sqrt{2 \beta^{-1} K} \dot{W}.
$$
The cell problem is
$$
\cL_0 \phi = \partial_y p
$$
and the homogenized coefficient is
$$
K =Z_p^{-1} \int_0^1 (1 + \partial_y \phi)^2  e^{-p(y)/\sigma} \, dy.
$$
We have that
$$
\langle p , q \rangle_{\alpha(x)} = \frac{\beta}{2 K} p q.
$$
The error in the likelihood is
$$
E_{\infty}(\theta) = \int_{-\infty}^{\infty} \int_0^1 \Big(\cL_1 \Gamma (x,y) - \langle
F, \cL_1 \phi \rangle_{\alpha(x)} \Big) \rho(x,y) \, dy dx,
$$
where
$$
\Gamma = \langle F, \phi \rangle_{\alpha(x)},
$$
$$
F = - K \partial_x V.
$$
We have that
$$
\Gamma(x,y) = \frac{\beta}{2 K} (-K \partial_x V \phi)
= - \frac{\beta}{2} \partial_x V \phi.
$$
Furthermore
$$
\cL_1 = -\partial_x V \partial_y - \partial_y p \partial_x +
2 \beta^{-1} \partial_x \partial_y.
$$
Consequently
\begin{eqnarray*}
\cL_1 \Gamma(x,y) = \frac{\beta}{2} \Big(|\partial_x V|^2 \partial_y \phi
+ \partial_y p \partial_x^2 V \phi - 2 \beta^{-1} \partial_x^2 V \partial_y \phi \Big).
\end{eqnarray*}
In addition,
$$
\langle F, \cL_1 \phi \rangle_{\alpha(x)} = \frac{\beta}{2} |\partial_x V |^2
\partial_y \phi.
$$
The error in the likelihood is
\begin{eqnarray*}
E_{\infty}(\theta) &=& \frac{\beta}{2} \int_{\R} \int_0^1 \Big( - \partial_y p
\partial_x^2 V \phi + 2 \beta^{-1} \partial_x^2 \partial_y \phi \Big) Z_V^{-1} Z_p^{-1}
e^{-\beta V(x;\theta) - \beta p(y)} \, dx dy \\ & = & - \frac{ Z_V^{-1} Z_p^{-1}}{2}
\int_{\R}
\partial_x^2 V e^{-\beta V(x;\theta)} \, dx \int_0^1 \partial_y\phi  e^{-\beta p(y)} \, dy
\\ &&+  Z_V^{-1} Z_p^{-1} \int_{\R} \partial_x^2 V e^{-\beta V(x;\theta)} \, dx \int_0^1
\partial_y \phi e^{-\beta p(y)} \, dy \\ & = & \frac{ Z_V^{-1} Z_p^{-1}}{2}
\int_{\R} \partial_x^2 V e^{-\beta V(x;\theta)} \, dx \int_0^1 \partial_y \phi e^{-\beta
p(y)} \, dy \\ & = & \frac{\beta Z_V^{-1}}{2} \int_{\R} |\partial_x V|^2 e^{-\beta
V(x;\theta)} \, dx \Big(-1 + \widehat{Z}_p^{-1} Z_p^{-1} \Big).
\end{eqnarray*}
In above derivation we used various integrations by parts, together with the formula for
the derivative of the solution of the Poisson equation $\partial_y \phi = -1 +
\widehat{Z}_p^{-1} e^{\beta p(y)}$, \cite[p. 213]{PavlSt08}. The fact that $E_{\infty}$
is nonpositive follows from the inequality $Z_p^{-1} \widehat{Z}_p^{-1} < 1$ (for $p(y)$
not identically equal to $0$), which follows from the Cauchy-Schwarz inequality. \qed
\begin{remark}
An application of Laplace's method shows that, for $\beta \gg 1$, $Z_p^{-1}
\widehat{Z}_p^{-1} \sim e^{- 2 \beta}$.
\end{remark}
\section{Conclusions}\label{sec:conc}
The problem of parameter estimation for fast/slow systems of SDEs which admit
a coarse-grained description in terms of an SDE for the slow variable was studied in this paper. It was shown that, when applied to the averaging
problem, the maximum likelihood estimator (MLE) is asymptotically unbiased
and we can use it to estimate accurately the parameters in the drift coefficient
of the coarse-grained
model using data from the slow variable in the fast/slow system. On the contrary, the MLE is asymptotically
biased when applied to the homogenization problem and a systematic asymptotic
error appears in the log-likelihood function, in the long time/infinite scale
separation limit. The MLE can lead to the correct estimation of the parameters
in the drift coefficient of the homogenized equation provided that we subsample the data from the fast/slow system  at the appropriate sampling rate.

The averaging/homogenization systems of SDEs that we consider in this paper are of quite
general form and have been studied quite extensively in the last several decades since
they appear in various applications, e.g. molecular dynamics, chemical kinetics,
mathematical finance, atmosphere/ocean science-see the references in~\cite{PavlSt08}.
Thus, we believe that our results show that great care has to be taken when using maximum
likelihood in order to infer information about parameters in stochastic systems with
multiple characteristic time scales.

There are various problems, both of theoretical and of applied interest,
that remain open and that we plan to address in future work. We list some
of them below.
\begin{itemize}
\item Bayesian techniques for parameter estimation of multiscale diffusion
processes.
\item The development of efficient algorithms for estimating the parameters
in the coarse-grained model of a fast/slow stochastic system. Based on the work that has
been done to similar models in the context of
econometrics~\cite{OlhPavlSyk08,AitMykZha05a} one expects that such an algorithm would
involve the estimation of an appropriate measure of scale separation $\eps$, and of the
optimal sampling rate, averaging over all the available data and a bias reduction step.
\item Investigate whether there is any advantage in using random sampling
rates.
\item Investigate similar issues for deterministic fast/slow systems of differential
equations.
\end{itemize}

\section*{Acknowledgements} AP has been partially supported by a  Marie Curie International Reintegration Grant, MIRG-CT-2005-029160. AMS is partially supported
by EPSRC.

\appendix
\section{Appendix}
\subsection{An Ergodic Theorem with Convergence Rates}
Consider the SDE
\begin{equation}
\frac{dz}{dt}=h(z)+\gamma(z)\frac{dW}{dt},
\label{eq:zeq}
\end{equation}
with $z \in \cZ$, where $\cZ$ is either $\bbR^k$ or $\bbT^k$, $h:\cZ \to \bbR^k,$
$\gamma:\cZ \to \bbR^{k \times p}$ and $w \in \bbR^p$ a standard Brownian motion.
Assume that $h, \gamma$ are $C^{\infty}$ with bounded derivatives. Let
$\psi: {\cal Z} \to \bbR$ be bounded, and $\phi:\cZ \to \bbR$ be bounded.
We denote the generator of the Markov process \eqref{eq:zeq} by ${\cal A}.$

\begin{assumptions} The equation \eqref{eq:zeq} is ergodic with invariant measure
$\nu(z)dz$. Let
$${\overline \phi}=\int_{\cal Z} \phi(z)\nu(z)dz.$$
Then the equation
$$-{\cal A}\Phi=\phi-{\overline \phi}, \quad \int_{\cZ}\Phi(z)\nu(z)dz=0$$
has a unique solution $\Phi:\cZ \to \bbR$, with $\Phi$ and $\nabla \Phi$ bounded.
\end{assumptions}

\begin{lemma} \label{l:A1} Let
$$I=\frac{1}{\sqrt T}\int_0^T \psi(z(t))dW(t).$$
Then there exists a constant $C>0$: $\bbE|I|^2 \le C$ for all $T>0.$
\end{lemma}

\begin{proof} Use the It\^{o} isometry and invoke the boundedness of $\psi.$
\end{proof}

\begin{lemma} \label{l:A2} Time averages converge to their mean value
almost surely. Furthermore  there is a constant $C>0$:
$$\bbE \left\vert\frac{1}{T}\int_0^T \phi(z(t))dt-{\overline \phi}\right\vert^2 \le
\frac{C}{T}.$$
\end{lemma}

\begin{proof}
By applying the It\^{o} formula to $\Phi$ we obtain
$$-\int_0^T {\cal A}\Phi(z(t))dt=\Phi(z(0))-\Phi(z(T))+\int_0^T
\left(\nabla\Phi \gamma \right)(z(t))dW(t).$$
Thus
\begin{align*}
\int_0^T \phi(z(t))dt&={\overline \phi}+\frac{1}{T}\left(\Phi(z(0))-\Phi(z(T))\right)+
\frac{1}{\sqrt T}I,\\
I&=\frac{1}{\sqrt T}\int_0^T \left(\nabla \Phi \gamma \right)(z(t))dW(t).
\end{align*}
The result concerning $L^2(\Omega)$ convergence follows from boundedness
of $\Phi$, $\nabla \Phi$ and $\gamma$, together with Lemma \ref{l:A1}.
Almost sure convergence follows from the ergodic theorem.
\end{proof}
%
\subsection{Consistency of the Estimators}

\begin{lemma}
\label{consistency}
Let $(\tilde{\Omega},\tilde{{\mathcal F}},\tilde{\mathbb P})$ be a probability space and $g_\epsilon:\tilde{\Omega}\times\Theta\rightarrow\R$, $g_0:\Theta\rightarrow\R$ be such that
\begin{equation}
\label{ass 1a}
\forall \theta\in\Theta,\ g_\epsilon\rightarrow g_0\ {\rm\ in\ probability,\ as}\ \epsilon\rightarrow 0
\end{equation}
and
\begin{equation}
\label{ass 1b}
\forall\delta,\kappa>0:\ {\mathbb P}\left\{\omega:\ \sup_{|u|>\delta}\left(g_\epsilon(\omega,\hat{\theta}_0 + u) - g_0(\hat{\theta}_0 + u)\right)>\kappa\ \right\}\rightarrow 0,\ {\rm as}\ \epsilon\rightarrow 0,
\end{equation}
where
\[ \hat{\theta}_0 = \arg\sup_{\theta\in\Theta} g_0(\theta).\]
Moreover, we assume that
\begin{equation}
\label{ass 2}
\forall\delta>0,\ \sup_{|u|>\delta}\left(g_0(\hat{\theta}_0+u)-g_0(\hat{\theta}_0)\right)\leq -\kappa(\delta)<0.
\end{equation}
If
\[ \hat{\theta}_\epsilon(\omega) = \arg\sup_{\theta\in\Theta}g_\epsilon(\omega,\theta)\]
then
\[ \hat{\theta}_\epsilon \rightarrow \hat{\theta}_0\ \ {\rm in\ probability}.\]
\end{lemma}
\begin{proof}
First note that $\forall\delta>0$
\begin{equation}
\label{inq 1}
\tilde{\mathbb P}\left\{ |\hat{\theta}_\epsilon-\hat{\theta}_0|>\delta\right\} \leq \tilde{\mathbb P}\left\{\ \sup_{|u|>\delta}\left(g_\epsilon(\omega,\hat{\theta}_0 + u) - g_\epsilon(\omega,\hat{\theta}_0)\right)\geq 0\ \right\}.
\end{equation}
We define
\[ G_\epsilon(\omega;\theta,u) := g_\epsilon(\omega,\theta + u) - g_\epsilon(\omega,\theta)\ {\rm and}\ G_0(\theta,u) := g_0(\theta + u) - g_0(\theta).\]
Clearly,
\[ \sup_{|u|>\delta} G_\epsilon(\omega;\hat{\theta}_0,u) \leq \sup_{|u|>\delta}\left( G_\epsilon(\omega;\hat{\theta}_0,u)-G_0(\hat{\theta}_0,u) \right) + \sup_{|u|>\delta}G_0(\hat{\theta}_0,u)\]
and thus
\begin{eqnarray}
\label{inq 2}
\nonumber \tilde{\mathbb P}\left\{ \sup_{|u|>\delta} G_\epsilon(\omega;\hat{\theta}_0,u)\geq 0\right\} &\leq&  \tilde{\mathbb P}\left\{\ \sup_{|u|>\delta}\left(G_\epsilon(\omega;\hat{\theta}_0,u) - G_0(\hat{\theta}_0,u)\right)\geq -\sup_{|u|>\delta}G_0(\hat{\theta}_0,u)\ \right\}\\
&\leq& \tilde{\mathbb P}\left\{\ \sup_{|u|>\delta}\left(G_\epsilon(\omega;\hat{\theta}_0,u) - G_0(\hat{\theta}_0,u)\right)\geq \kappa(\delta)>0\ \right\}
\end{eqnarray}
by Assumption (\ref{ass 2}). Note that
\[ G_\epsilon(\omega;\hat{\theta}_0,u) - G_0(\hat{\theta}_0,u) = \left( g_\epsilon(\omega;\hat{\theta}_0+u) - g_0(\hat{\theta}_0+u) \right) - \left( g_\epsilon(\omega;\hat{\theta}_0) - g_0(\hat{\theta}_0) \right).\]
So, by conditioning on $\left\{ \omega:\ |g_\epsilon(\omega;\hat{\theta}_0) - g_0(\hat{\theta}_0)|\geq \frac{1}{2}\kappa(\delta)\right\}$ and (\ref{inq 1}) and (\ref{inq 2}), we get that
\begin{eqnarray*}
\tilde{\mathbb P}\left\{ |\hat{\theta}_\epsilon-\hat{\theta}_0|>\delta\right\} \leq&
\tilde{\mathbb P}\left\{\ \sup_{|u|>\delta}\left(g_\epsilon(\omega;\hat{\theta}_0+u) - g_0(\hat{\theta}_0+u)\right)\geq \frac{1}{2}\kappa(\delta)>0\ \right\} \\
&+
\tilde{\mathbb P}\left\{\ |g_\epsilon(\omega;\hat{\theta}_0) - g_0(\hat{\theta}_0)|\geq \frac{1}{2}\kappa(\delta)>0\right\}
\end{eqnarray*}
Both probabilities on the right-hand-side go to zero as $\epsilon\rightarrow 0$, by assumptions (\ref{ass 1b}) and (\ref{ass 1a}) respectively. We conclude that $\hat{\theta}_\epsilon \rightarrow \hat{\theta}_0$ in probability.
\end{proof}
%
%
%
\subsection{Proof of Propositions~\ref{cor:x_estim} and~\ref{lem:ergodic}}
In this section we present the proofs of Propositions~\ref{cor:x_estim} and
\ref{lem:ergodic} which we repeat there, for the reader's convenience.
\begin{prop}
\label{A:5}
Let $(x(t),y(t))$ be the solution of~\eqref{eq:hom} and assume that
Assumptions~\ref{ass:a1} and~\ref{ass:a2} hold. Then, for $\eps, \, \delta$ sufficiently
small, the increment of the process $x(t)$ can be written in the form
\begin{equation*}
x_{n+1} - x_n = F(x_n;\theta_0) \, \delta + M_n + R(\eps, \delta),
\end{equation*}
where $M_n$ denotes the martingale term
$$
M_n  =  \int_{n \delta}^{(n+1) \delta} \left(\nabla_y \Phi\beta + \alpha_0 \right) (x(s),
y(s)) \, d V(s) + \int_{n \delta}^{(n+1) \delta} \alpha_1 (x (s) , y(s)) \, d U(s)
$$
with $\|M_n \|_p \leq C \sqrt{\delta}$ and
$$
\| R(\eps, \delta) \|_p \leq C (\delta^{3/2} + \eps \delta^{\frac12} + \eps).
$$
\end{prop}

\smallskip

\begin{prop}
\label{A:6}
Let $g \in C^1(\cX)$ and let Assumptions \ref{ass:a3} hold. Assume that $\eps$ and $N$
are related as in Theorem \ref{th:main1}. Then
\begin{equation*}
\lim_{\eps \to 0} \frac{1}{N}\sum_{n=0}^{N-1} g(x_n)=\bbE^{\pi} g,
\end{equation*}
where the convergence is in $L^2$ with respect to the measure on initial conditions with
density $\pi(x)\rho(y;x).$
\end{prop}

For the proofs of Propositions~\ref{A:5} and~\ref{A:6}, both used in the proof of
Theorem~\ref{th:main1}, we will need the following two technical lemmas. We start with a
rough estimate on the increments of the process $x (t)$.
\begin{lemma}\label{lem:estim_increm}
Let $(x(t),y(t))$ be the solution of~\eqref{eq:hom} and assume that
Assumptions~\ref{ass:a1} and~\ref{ass:a2} hold. Let $s \in [n\delta,(n+1)\delta].$ Then,
for $\eps, \, \delta$ sufficiently small, the following estimate holds:
\begin{equation}\label{e:estim_increm}
\|x(s) - x_n \|_p \leq C (\eps + \delta^{\frac12}).
\end{equation}
\end{lemma}
\proof We apply It\^{o}'s formula to $\Phi$, the solution of the Poisson
equation~\eqref{e:poisson_main}, to obtain
\begin{eqnarray*}
x(s) -x_n & = & -\eps (\Phi(x(s), y(s)) - \Phi(x_{n}, y_{n})) + \int_{n \delta}^{s}
\left( \cL_1 \Phi + f_1) \right)(x (s) , y(s)) \, ds
\\ &&
+ \int_{n \delta}^{s} \left( \nabla_y \Phi \beta+ \alpha_0 \right) (x (s) , y(s)) \, d
V(s) + \int_{n \delta}^{s} \alpha_1 (x (s) , y(s)) \, d U(s)
\\ &&
+ \eps \int_{n \delta}^{s} ( \cL_2 \Phi )(x (s) , y(s)) \, ds + \eps \int_{n \delta}^{s}
\left( \nabla_y \Phi \alpha_0\right)(x (s) , y(s)) \, dU(s)
\\ &&
+ \eps \int_{n \delta}^{s} \left(\nabla_x \Phi \alpha_1\right) (x (s) , y(s)) \, d V(s)
\\ & =: & J_1 + J_2 + J_3 + J_4 + J_5 + J_6 + J_7.
\end{eqnarray*}
Our assumptions on $\Phi(x,y)$, together with standard inequalities, imply that
\begin{eqnarray*}
&& \|J_1 \|_p \leq C \eps, \; \; \|J_2 \|_p \leq C \delta, \; \; \|J_3 \|_p \leq C
\delta^{\frac12},
\\ &&
\|J_4 \|_p \leq C \delta^{\frac12}, \; \; \|J_5 \|_p \leq C \eps \delta, \; \; \|J_6 \|_p
\leq C \eps \delta^{1/2}, \; \; \|J_7 \|_p \leq C \eps \delta^{1/2}.
\end{eqnarray*}
Estimate~\eqref{e:estim_increm} follows from these estimates. \qed

\smallskip

Using this lemma we can prove the following estimate.

\begin{lemma}\label{prop:estim_int}
Let $h(x,y)$ be a smooth, bounded function, let $(x(t),y(t))$ be the solution
of~\eqref{eq:hom} and assume that Assumption~\ref{ass:a1} holds. Define
$$
H(x) := \int_{\cY} h(x,y) \, \rho(y;x)dy.
$$
Then, for $\eps, \, \delta$ sufficiently small, the
following estimate holds:
\begin{equation}\label{e:estim_int}
\int_{n \delta}^{(n+1) \delta} h (x (s) , y (s)) \, ds = H(x_n) \, \delta + R(\eps,
\delta)
\end{equation}
where
$$
\| R(\eps, \delta)  \|_p \leq C (\eps^2 + \delta^{3/2} + \eps \delta^{1/2}).
$$
\end{lemma}

\proof Let $\phi$ be the mean zero solution of the
equation
\begin{equation}\label{e:poisson}
- \cL_0 \phi = h(x,y) - H(x).
\end{equation}
By Assumption~\ref{ass:a1} this solution is smooth in both $x, \, y$ and it is unique and
bounded. We apply It\^{o}'s formula to obtain
\begin{eqnarray*}
\int_{n \delta}^{(n+1) \delta} \left( h(x(s) , y(s))) - H(x (s)) \right) \, ds
& = & -\eps^2 (\phi(x_{n+1}, y_{n+1})  - \phi(x_{n}, y_{n}))\\
&& + \eps \int_{n \delta}^{(n+1) \delta} \cL_1 \phi(x(s) , y (s)) \, ds \\ && + \eps^2
\int_{n \delta}^{(n+1) \delta}  \cL_2 \phi(x (s) , y(s)) \, ds \\&&+ \eps^2 \int_{n
\delta}^{(n+1) \delta} (\nabla_x \phi\alpha_0) (x
(s) , y(s)) \, d U(s)  \\
&& + \eps \int_{n \delta}^{(n+1) \delta} (\nabla_y \phi\beta + \eps \nabla_x \phi
\alpha_1)(x (s) , y(s)) \, d V(s) \\ & =: & J_1 + J_2 + J_3 + J_4 + J_5.
\end{eqnarray*}
Our assumptions on the solution $\phi$ of the Poisson equation~\eqref{e:poisson},
together with standard estimates for the moments of stochastic integrals and H\"{o}lder's
inequality give the estimates
\begin{eqnarray*}
&& \|J_1 \|_p \leq C \eps^2, \; \; \|J_2 \|_p \leq C \eps \delta, \; \; \|J_3 \|_p \leq C
\eps^2 \delta,
\\ &&
\|J_4 \|_p \leq C \eps^2 \delta^{1/2}, \; \; \|J_5 \|_p \leq C \eps
\delta^{1/2}.
\end{eqnarray*}
The above estimates imply that
\begin{eqnarray*}
\int_{n \delta}^{(n+1) \delta} h(x (s) , y(s)) \, ds = \int_{n \delta}^{(n+1) \delta} H(x
(s)) \, ds + R_1(\eps, \delta)
\end{eqnarray*}
with
$$
\| R_1(\eps, \delta) \|_p \leq C \left( \eps \delta^{1/2} + \eps^2 \right).
$$
We use the H\"{o}lder inequality and the Lipschitz continuity of $H(x)$
to estimate:
\begin{eqnarray*}
\left\| \int_{n \delta}^{(n+1) \delta} H(x (s)) \, ds  - H(x_n) \, \delta \right\|_p^p &
= & \left\| \int_{n \delta}^{(n+1) \delta} \left( H(x (s))  - H(x_n) \right) \, ds
\right\|_p^p \\ & \leq & \delta^{p-1}\int_{n \delta}^{(n+1) \delta} \left\| H(x (s)) -
H(x_n) \right\|_p^p \, ds \\ & \leq & C \delta^{p-1}\int_{n \delta}^{(n+1) \delta}
\left\| x (s)  - x_n \right\|_p^p \, ds
\\ & \leq & C\delta^{p} \left( \delta^{1/2} + \eps \right)^p = R_2(\eps, \delta)^p,
\end{eqnarray*}
where Lemma~\ref{lem:estim_increm} was used and $R_2(\eps,\delta)=(\eps \delta
+\delta^{3/2})$.
We combine the above estimates to obtain
\begin{eqnarray*}
\int_{n \delta}^{(n+1) \delta} h(x(s) , y(s))) \, ds & = & \int_{n \delta}^{(n+1) \delta}
H(x (s)) \, ds + R_1(\eps, \delta) \\ & = & H(x_n) \, \delta + R_1(\eps, \delta) +
R_2(\eps, \delta),
\end{eqnarray*}
from which~\eqref{e:estim_int} follows. \qed

\medskip

{\it Proof of Proposition~\ref{cor:x_estim} (Proposition \ref{A:5}).}
This follows from the first line of the proof of Lemma
\ref{lem:estim_increm}, the estimates therein concerning all the $J_i$ with
the exception of $J_2$, and the use of Lemma \ref{prop:estim_int} to
estimate $J_2$ in terms of $\delta F(x_n;\theta_0).$
\qed

\medskip

{\it Proof of Proposition~\ref{lem:ergodic} (Proposition \ref{A:6}).} We have
\begin{eqnarray*}
\frac{1}{N} \sum_{n=0}^{N-1} g(x_n) & = & \frac{1}{N \delta} \sum_{n=0}^{N-1} \int_{n
\delta}^{(n+1) \delta}g(x_n) \, ds \\ & = & \frac{1}{N \delta} \sum_{n=0}^{N-1} \int_{n
\delta}^{(n+1) \delta}g(x(s)) \, ds + \frac{1}{N \delta} \sum_{n=0}^{N-1} \int_{n
\delta}^{(n+1)\delta} \left( g(x_n) - g(x(s)) \right)
\, ds\\
&=& \frac{1}{N \delta} \int_0^{N\delta} g(x(s))ds+\frac{1}{N \delta} \sum_{n=0}^{N-1}
\int_{n \delta}^{(n+1)\delta} \left( g(x_n) - g(x(s)) \right)
\, ds\\
&=:& I_1+ R_1.
\end{eqnarray*}

We introduce the notation
$$
f_n := \int_{n \delta}^{(n+1) \delta} \big(g(x_n) - g(x(s)) \big) \, ds.
$$
By Lemma~\ref{lem:estim_increm} we have that $x(s)-x_n={\cal O}(\eps+\delta^{\frac12})$
in $L^p(\Omega').$ We use this, together with the Lipschitz continuity of $g$ and
H\"{o}lder's inequality, to estimate:
\begin{eqnarray*}
\|f_n \|_p^p & \leq & \delta^{p/q}
\int_{n \delta}^{(n+1) \delta} \E \big|g(x_n) - g(x(s)) \big|^p \, ds \\
& \leq & C \delta^{1 + p/q} \big( \eps^p + \delta^{p/2} \big).
\end{eqnarray*}
Here $p^{-1}+q^{-1}=1.$ Using this we can estimate $R_1$ using:
\begin{eqnarray*}
\| R_1 \|_p & \leq & \frac{1}{N \delta} \sum_{n=0}^{N-1} \|f_n \|_p \leq C \frac{1}{N
\delta} N \delta^{(1/p + 1/q)} \big( \eps + \delta^{1/2} \big)=
C \big( \eps + \delta^{1/2} \big)\rightarrow 0,
\end{eqnarray*}
as $\eps \rightarrow 0$.

Thus it remains to estimate $I_1$. Let $T = N \delta$.  Let $\Phie$ solve
\begin{equation}\label{e:phie}
-{\cal L}_{hom}\Phie(x,y)=\hat{g}(x):=g(x)-\bbE^{\pie}g.
\end{equation}
Apply It\^{o}'s formula. This gives
\begin{align*}
\frac{1}{T}\int_0^T g(x(s))ds-\bbE^{\pie}g=& -\frac{1}{T}\Bigl(\Phie \big( x(T),y(T)
\big)
-\Phie \big( x(0),y(0)\big)  \Bigr)\\
&+\frac{1}{\eps T}\int_0^T \bigl(\nabla_y\Phie \beta)(x(s),y(s)) \,
dV(s)\\
&+\frac{1}{T}\int_0^T \bigl(\nabla_x\Phie \alpha)(x(s),y(s))
\, dU'(s),\\
& =: J_1+J_2
\end{align*}
where $J_2$ denotes the two stochastic integrals and we write $\alpha dU'=\alpha_0
dU+\alpha_1 dV,$ in law. Note that
$$\bbE^{\pie}g \to \bbE^{\pi} g$$
as $\eps \to 0$ by Assumptions \ref{ass:a3}.
Thus the theorem will be proved if we can show
that $J_1+J_2$ tends to zero in the required topology on the initial
conditions. Note that
\begin{align*}
\bbE^{\rho^{\eps}}|J_1|^2 & \le \frac{4}{T^2} \bbE^{\rho^{\eps}}|\psi^{\eps}|^2,\\
\bbE^{\rho^{\eps}}|J_2|^2 & \le \frac{1}{T} \bbE^{\rho^{\eps}}
\langle \nabla \psi^{\eps},\Sigma\nabla \psi^{\eps} \rangle.
\end{align*}
Here $\Sigma$ is defined in Assumptions \ref{assump:a4} and $\nabla$ is the gradient with
respect to $(x^T,y^T)^T.$ We note that, by stationarity, we have that
\begin{equation}\label{e:notation}
\E^{\rho^{\eps}} |\psi^{\eps}|^2 = \| \psi^{\eps}\|, \quad \bbE^{\rho^{\eps}} \langle \nabla
\psi^{\eps},\Sigma\nabla \psi^{\eps} \rangle = (\nabla \psi^{\eps},\Sigma\nabla
\psi^{\eps}),
\end{equation}
where $\|\cdot \|$ and $( \cdot, \cdot )$ denote the $L^2(\cX \times \cY; \mu^{\eps}(dx
dy))$ norm and inner product, respectively.

Use of the Dirichlet form (see Theorem 6.12 in \cite{PavlSt08}) shows that
\begin{align*}
\big( \nabla \psi^{\eps},\Sigma\nabla \psi^{\eps} \big)
& \le 2\int {\hat g}(x)\psi^{\eps}(x,y)\rho^{\eps}(x,y)dxdy\\
& \le a \|{\hat g}\|^2+a^{-1}  \|\psi^{\eps}\|^2,
\end{align*}
for any $a>0.$
Using the Poincar\'{e} inequality \eqref{e:poincare}, together with Assumptions
\ref{ass:a3} and \ref{assump:a4}, gives
$$
\|\psi^{\eps}\|^2 \le C_p^2 \|\nabla \psi^{\eps}\|^2 \le  a C_{\gamma}^{-1} C_p^2
\|{\hat g}\|^2+a^{-1} C_{\gamma}^{-1} C_p^2 \|\psi^{\eps}\|^2.
$$
Choosing $\alpha$ so that $a^{-1} C_{\gamma}^{-1} C_p^2= \frac{1}{2}$ gives
$$
\|\psi^{\eps}\|^2 \le C \bbE^{\rho^{\eps}} |{\hat g}|^2.
$$
Hence
$$
 \big( \psi^{\eps},\Gamma \nabla \psi^{\eps} \big) \le C \bbE^{\rho^{\eps}} |{\hat g}|^2,
$$
where the notation introduced in~\eqref{e:notation} was used. The constant $C$ in the
above inequalities is independent of $\eps$. Thus
\begin{equation}
\bbE^{\rho^{\eps}}|J_1|^2 +\bbE^{\rho^{\eps}}|J_2|^2  \le \frac{1}{T} C
\bbE^{\rho^{\eps}} |{\hat g}|^2.
\end{equation}
Since the measure with density $\rho^{\eps}$ converges to the measure with
density $\pi(x)\rho(y;x)$ the desired result follows. \qed
%

\def\cprime{$'$} \def\cprime{$'$} \def\cprime{$'$} \def\cprime{$'$}
  \def\cprime{$'$} \def\cprime{$'$} \def\cprime{$'$}
  \def\Rom#1{\uppercase\expandafter{\romannumeral #1}}\def\u#1{{\accent"15
  #1}}\def\Rom#1{\uppercase\expandafter{\romannumeral #1}}\def\u#1{{\accent"15
  #1}}\def\cprime{$'$} \def\cprime{$'$} \def\cprime{$'$} \def\cprime{$'$}
  \def\cprime{$'$} \def\cprime{$'$} \def\cprime{$'$}


\begin{thebibliography}{10}

\bibitem{AitMykZha05b}
Y.~Ait-Sahalia, P.~A. Mykland, and L~Zhang.
\newblock How often to sample a continuous-time process in the presence of
  market microstructure noise.
\newblock {\em Rev. Financ. Studies}, 18:351--416, 2005.

\bibitem{AitMykZha05a}
Y.~Ait-Sahalia, P.~A. Mykland, and L~Zhang.
\newblock A tale of two time scales: Determining integrated volatility with
  noisy high-frequency data.
\newblock {\em J. Amer. Stat. Assoc.}, 100:1394--1411, 2005.

\bibitem{BakryCattiauxGuillin08}
D.~Bakry, P.~Cattiaux, and A.~Guillin.
\newblock Rate of convergence for ergodic continuous {M}arkov processes:
  {L}yapunov versus {P}oincar\'e.
\newblock {\em J. Funct. Anal.}, 254(3):727--759, 2008.

\bibitem{lions}
A.~Bensoussan, J.-L. Lions, and G.~Papanicolaou.
\newblock {\em Asymptotic analysis for periodic structures}, volume~5 of {\em
  Studies in Mathematics and its Applications}.
\newblock North-Holland Publishing Co., Amsterdam, 1978.

\bibitem{dudley84}
R.~M. Dudley.
\newblock A course on empirical processes.
\newblock In {\em \'Ecole d'\'et\'e de probabilit\'es de Saint-Flour,
  XII---1982}, volume 1097 of {\em Lecture Notes in Math.}, pages 1--142.
  Springer, Berlin, 1984.

\bibitem{ELV05}
W.~E, D.~Liu, and E.~Vanden-Eijnden.
\newblock Analysis of multiscale methods for stochastic differential equations.
\newblock {\em Comm. Pure Appl. Math.}, 58(11):1544--1585, 2005.

\bibitem{EthKur86}
S.N. Ethier and T.G. Kurtz.
\newblock {\em Markov processes}.
\newblock Wiley Series in Probability and Mathematical Statistics: Probability
  and Mathematical Statistics. John Wiley \& Sons Inc., New York, 1986.

\bibitem{GivKevKup06}
D.~Givon, I.G. Kevrekidis, and R.~Kupferman.
\newblock Strong convergence schemes of projective intregration schemes for
  singularly perturbed stochastic differential equations.
\newblock {\em Comm. Math. Sci.}, 4(4):707--729, 2006.

\bibitem{GKS04}
D.~Givon, R.~Kupferman, and A.M. Stuart.
\newblock Extracting macroscopic dynamics: model problems and algorithms.
\newblock {\em Nonlinearity}, 17(6):R55--R127, 2004.

\bibitem{KSh91}
I.~Karatzas and S.E. Shreve.
\newblock {\em Brownian {M}otion and {S}tochastic {C}alculus}, volume 113 of
  {\em Graduate Texts in Mathematics}.
\newblock Springer-Verlag, New York, second edition, 1991.

\bibitem{Kut04}
Y.A. Kutoyants.
\newblock {\em Statistical inference for ergodic diffusion processes}.
\newblock Springer Series in Statistics. Springer-Verlag London Ltd., London,
  2004.

\bibitem{nelson}
E.~Nelson.
\newblock {\em Dynamical theories of {B}rownian motion}.
\newblock Princeton University Press, Princeton, N.J., 1967.

\bibitem{OlhPavlSyk08}
S.~Olhede, G.A. Pavliotis, and A.~Sykulski.
\newblock Multiscale inference for high frequency data.
\newblock {\em Preprint}, 2008.

\bibitem{PapStrVar77}
G.~C. Papanicolaou, D.W. Stroock, and S.~R.~S. Varadhan.
\newblock Martingale approach to some limit theorems.
\newblock In {\em Papers from the Duke Turbulence Conference (Duke Univ.,
  Durham, N.C., 1976), Paper No. 6}, pages ii+120 pp. Duke Univ. Math. Ser.,
  Vol. III. Duke Univ., Durham, N.C., 1977.

\bibitem{PV01}
E.~Pardoux and A.~Yu. Veretennikov.
\newblock On the {P}oisson equation and diffusion approximation. {I}.
\newblock {\em Ann. Probab.}, 29(3):1061--1085, 2001.

\bibitem{PV03}
{\`E}.~Pardoux and A.~Yu. Veretennikov.
\newblock On {P}oisson equation and diffusion approximation. {II}.
\newblock {\em Ann. Probab.}, 31(3):1166--1192, 2003.

\bibitem{PV05}
E.~Pardoux and A.~Yu. Veretennikov.
\newblock On the {P}oisson equation and diffusion approximation. {III}.
\newblock {\em Ann. Probab.}, 33(3):1111--1133, 2005.

\bibitem{PavlSt03}
G.~A. Pavliotis and A.~M. Stuart.
\newblock White noise limits for inertial particles in a random field.
\newblock {\em Multiscale Model. Simul.}, 1(4):527--533 (electronic), 2003.

\bibitem{PavlSt06}
G.~A. Pavliotis and A.~M. Stuart.
\newblock Parameter estimation for multiscale diffusions.
\newblock {\em J. Stat. Phys.}, 127(4):741--781, 2007.

\bibitem{PavlSt08}
G.A. Pavliotis and A.M. Stuart.
\newblock {\em Multiscale methods}, volume~53 of {\em Texts in Applied
  Mathematics}.
\newblock Springer, New York, 2008.
\newblock Averaging and homogenization.

\bibitem{Rao99}
B.~L.~S. Prakasa~Rao.
\newblock {\em Statistical inference for diffusion type processes}, volume~8 of
  {\em Kendall's Library of Statistics}.
\newblock Edward Arnold, London, 1999.

\bibitem{vanZanten01}
J.~H. van Zanten.
\newblock A note on consistent estimation of multivariate parameters in ergodic
  diffusion models.
\newblock {\em Scand. J. Statist.}, 28(4):617--623, 2001.

\bibitem{Vand03}
E.~Vanden-Eijnden.
\newblock Numerical techniques for multi-scale dynamical systems with
  stochastic effects.
\newblock {\em Commun. Math. Sci.}, 1(2):385--391, 2003.

\bibitem{Vil04HPI}
C.~Villani.
\newblock {\em Hypocoercivity}.
\newblock AMS, 2008.

\end{thebibliography}
\end{document}